\documentclass[12 pt]{article}
\usepackage{fullpage}
\usepackage{lipsum}
\usepackage{amsfonts}
\usepackage{graphicx}
\usepackage{epstopdf}
\usepackage{standalone}
\usepackage{algcompatible}

\ifpdf
\DeclareGraphicsExtensions{.eps,.pdf,.png,.jpg}
\else
\DeclareGraphicsExtensions{.eps}
\fi

\usepackage{amsopn}

\usepackage{epsfig} 
\usepackage{amsmath} 
\usepackage{amssymb}  
\usepackage{bm}
\DeclareMathAlphabet{\mathcal}{OMS}{cmsy}{m}{n}
\usepackage[english]{babel}
\usepackage{url}
\usepackage{algorithm}
\usepackage{color}
\usepackage{pgfplots}
\usepackage{siunitx}
\usepackage{tikz}
\usepackage{tkz-euclide}
\usepackage{chngcntr}
\usepackage{verbatim}
\usepackage{enumerate}

\definecolor{ao(english)}{rgb}{0.0, 0.5, 0.0}
\usepackage[colorlinks, citecolor = {ao(english)}, linkcolor = {ao(english)}]{hyperref} 
\usepackage{cleveref}
\usepackage{aliascnt}
\usetikzlibrary{calc}
\usepackage{subcaption}
\usepackage{tikz}
\usepackage{pgfplots}
\usetikzlibrary{arrows,shapes,trees,calc,positioning,patterns,decorations.pathmorphing,decorations.markings}
\usetikzlibrary{matrix}
\usepgfplotslibrary{groupplots}
\pgfplotsset{compat=newest}
\usepackage{amsthm}
\usepgfplotslibrary{patchplots}

\crefname{figure}{Fig.}{Fig.}

\newtheorem{thm}{Theorem}
\crefname{thm}{Theorem}{Theorems}

\newtheorem{prop}{Proposition}
\crefname{prop}{Proposition}{Propositions}
\newtheorem{lem}{Lemma}
\crefname{lem}{Lemma}{Lemmas}
\newtheorem{cor}{Corollary}
\crefname{cor}{Corollary}{Corollaries}
\theoremstyle{remark}

\crefname{rem}{Remark}{Remarks}

\theoremstyle{definition}

\crefname{example}{Example}{Examples}

\crefname{ass}{Assumption}{Assumption}
\usepackage{dsfont}
\let\mathbb=\mathds
\usepackage[numbers,sort&compress]{natbib}

\crefname{conj}{Conjecture}{Conjectures}

\theoremstyle{definition}
\newtheorem{defn}{Definition}
\crefname{defn}{Definition}{Definitions}

\crefname{prob}{Problem}{Problems}
\crefname{algorithm}{Algorithm}{Algorithms}

\usepackage{dsfont}
\let\mathbb=\mathds
\usepackage[numbers,sort&compress]{natbib}

\newcommand{\Rmn}{\mathbb{R}^{n \times m}}

\newcommand{\Rnn}{\mathbb{R}^{n \times n}}

\newcommand{\rk}{\textnormal{rank}}

\newcommand{\diag}{\textnormal{diag}}
\newcommand{\sign}{\textnormal{sign}}

\newcommand{\argmin}{\operatornamewithlimits{argmin}}

\newcommand{\transp}{\mathsf{T}}
\newcommand{\vari}[1]{\text{S}^{-}(#1)}

\newcommand{\variz}[1]{S^{+}(#1)}
\newcommand{\sr}[1]{{SR}_{#1}}
\newcommand{\sgc}[1]{{SC}_{#1}}

\newcommand{\ssgc}[1]{{SSC}_{#1}}
\newcommand{\ovd}[1]{\text{OVD}_{#1}}

\newcommand{\tp}[1]{{#1}\text{-positive}}

\newcommand{\Con}[1]{{\mathcal{C}^{#1}}}
\newcommand{\Obs}[1]{{\mathcal{O}^{#1}}}

\newcommand{\compset}[1]{{#1}^\text{c}}

\newcommand{\vd}[1]{\text{VD}_{#1}}
\newcommand{\vb}[1]{\text{VB}_{#1}}
\newcommand{\svb}[1]{\text{SVB}_{#1}}

\newcommand{\nset}[1]{(1:#1)}
\newcommand{\submatrix}[3]{{#1}_{#2,#3}}

\newcommand{\compound}[2]{#1_{[#2]}}
\newcommand{\pcompound}[2]{{\left(#1\right)}_{[#2]}}

\colorlet{FigColor1}{blue}
\colorlet{FigColor2}{red}
\colorlet{FigColor3}{ao(english)}
\colorlet{FigColor4}{orange}
\pgfplotsset{every axis plot/.append style={line width=1.5pt}}
\begin{document}

\title{On System Operators with Variation Bounding Properties}
\author{ Christian Grussler\thanks{The author is with the Stephen B. Klein Faculty of Aerospace Engineering at the Technion -- Israel Institute of Technology, 3200003 Haifa, Israel.
		{\tt\small cgrussler@technion.ac.il}} \quad Chaim Roth\thanks{The author was with the faculty of Mechanical Engineering at the Technion -- Israel Institute of Technology, 3200003 Haifa, Israel. {\tt\small chaim.roth@campus.technion.ac.il}} \quad Kang Tong\thanks{The author is with the Stephen B. Klein Faculty of Aerospace Engineering at the Technion -- Israel Institute of Technology, 3200003 Haifa, Israel.
		{\tt\small kang.tong@technion.ac.il}}}%

\maketitle

\begin{abstract}
The property of linear discrete-time time-invariant system operators mapping inputs with at most $k-1$ sign changes to outputs with {at most} $k-1$ sign changes is investigated. We show that this property is tractable via the notion of $k$-sign consistency in case of the observability/controllability operator, which as such can also be used as a sufficient condition for the Hankel operator. Our results complement the mathematical literature by providing an algebraic characterization, independent of rank and dimension for variation bounding and diminishing matrices as well as by discussing their computational tractability. Based on these, we conduct our studies of variation bounding system operators beyond existing studies on order-preserving $k$-variation diminishment. Our findings are applied to the open problem of bounding the number of sign changes in a system's impulse response, {which appears, e.g., when bounding the number of over- and undershoots in a step response or the number of bangs in bounded optimal control problems. }
\end{abstract}

\section{Introduction}
Linear time-invariant (LTI) systems {$(A,b,c)$ given by}
\begin{equation}\label{eq:SISO_d}
    \begin{aligned}
 	    x(t+1) &= A x(t) + b u(t) \\
        y(t) & = cx(t),
    \end{aligned}
\end{equation}
$A \in \mathds{R}^{n \times n}$, $b,c^\transp \in \mathbb{R}^n$, that map nonnegative inputs $u$ to nonnegative outputs $y$ are characterized by \emph{nonnegative impulse} responses $g(t) := cA^{t-1}b \geq 0$, $t \geq 1$ and referred to as \emph{externally positive} \cite{farina2011positive}. In the particular case that $A$, $b$ and $c$ are element-wise nonnegative, the system is called \emph{internally positive} as it is externally positive and preserves the positivity of the state for nonnegative inputs. Systems of this kind frequently occur through compartmental modelling \cite{brown1980compartmental,farina2011positive} and are known for several advantageous analytical properties such as scalable stability certificates \cite{rantzer2015scalable}, optimal controller design \cite{tanaka2011bounded,pates2019optimal} or the  avoidance of over- and undershooting in closed-loop design  \cite{grussler2019tractable,darbha2003synthesis,phillips1988conditions,taghavian2023external}.  

In recent years, efforts have been made to study systems with extended positivity properties based on their influence on the variation, i.e., number of sign changes, of the input signal $u$ or the state $x(t)$ (see, e.g., \cite{weiss2019generalization,margaliot2018revisiting,grussler2020variation,grussler2021internally}. A mapping $u \mapsto Xu$ is called \emph{$k$-variation bounding} ($\vb{k}$) if it maps an input $u$ with at most $k$ sign changes to an output $Xu$ of at most $k$ sign changes. If the mapping is $\vb{j}$ for all $j \in \{0,\dots,k\}$, it is called \emph{$k$-variation diminishing} ($\vd{k}$) and if, additionally, the order of sign changes is preserved whenever $u$ and $Xu$ have equal variation, the $\vd{k}$ property is called \emph{order-preserving} ($\ovd{k}$). 

Systems with input-output $\ovd{k}$ properties have been studied in terms of Hankel and Toeplitz operator \cite{grussler2020variation}, which resulted in intermediates between external positivity ($\ovd{0}$ operators) and the well-known class of \emph{relaxation systems} \cite{willems1976realization} ($\ovd{\infty}$ Hankel operators). These intermediates have been exploited in \cite{grussler2020balanced} to provide guarantees for internal positivity preserving model order reduction -- a property well-known for $\ovd{\infty}$ Hankel operators. Studies of autonomous (unforced) systems with $\vd{k}$ and $\vb{k}$ state-transition maps \cite{margaliot2018revisiting,weiss2021generalization,ofir2024contraction} have provided an extended, unified approach to several important non-linear system analysis tools such as monotonicity/cooperativity, contraction analysis or Poincar{\'e}-Bendixson results. A system class that was proposed to merge the input-output and autonomous case has been suggested in \cite{grussler2021internally}. This resulted in a first certificate for $\ovd{k}$ observability/controllability operators and its application to bounding the variation of an impulse response beyond external positivity. {While the latter problem is of interest, e.g., when bounding the number of over- and undershoots in a step response or the number of bangs in bounded optimal control problems \cite{luenberger1968optimization,friesz2010dynamic}, it is still an open problem, where several lower bounds \cite{damm2014zero,swaroop1996some,elkhoury1993discrete} have been derived, but only few upper bounds \cite{elkhoury1993discrete}. }

The goal of this work is to derive a tractable characterization of systems with $\vb{k}$ observability/controllability operators, which complements the literature in several aspects: firstly, we close the gap of algebraic characterization of $\vb{k}$ and $\vd{k}$ matrices beyond limitations to rank and dimensions (cf. \cite{weiss2019generalization,karlin1968total}) or $\ovd{k}$ matrices (cf. \cite{karlin1968total,grussler2020variation}). Based on that, computationally tractable certificates are achieved by extending related results in \cite{pena_matrices_1995}. Secondly, it is shown that an application of these results to the infinite-dimensional observability/controllability operator leads to the external positivity of a related family of linear time-invariant systems. This is similar to the approach in \cite{grussler2020variation}, but goes beyond $\ovd{k}$ and does not require $A$ to be $\ovd{k}$ (cf. \cite{grussler2021internally}). To underline this benefit, we present a system example for which there exists a realization with $\vd{2}$ observability operator, but no realization with an $\vd{2}$ $A$. This parallels the fact that not all externally positive systems have an internally positive realization \cite{benvenuti2004tutorial,farina2011positive,grussler2019tractable}. Further, as any Hankel operator factors into observability and controllability operator, our results also provide a first sufficient certificate of $\vb{k}$ and $\vd{k}$ Hankel operators. Our results are illustrated by bounding the variation of the impulse response aka. the number of over- and undershoots in a step response. In future work, we hope to use these insights to also derive characteristics in terms of zeros and poles. {Lastly, the presented results have been instrumental to subsequent works  dealing with sparse optimization problems \cite{marmary2025tractable} and self-oscillations in relay feedback systems \cite{grussler2025discrete,tong2025selfsustained}.}

Finally note that this paper is an extension of the authors' earlier work \cite{roth2024variation}. The present work differs from \cite{roth2024variation} in several major aspects: (i) previously only the special case of strict
variation bounding was considered, which could be built upon existing matrix theory results \cite{pena_matrices_1995,karlin1968total}. As mentioned earlier, the current manuscript extends the system operator and matrix theory results substantially to the non-strict case. In particular, much weaker assumptions on rank and dimensions are made compared to \cite{weiss2019generalization,karlin1968total}. These novel results are valuable independent of the considered application and also directly include the case of (non-order-preserving) $k$-variation diminishment; (ii) the current paper includes an additional example to illustrate further benefits and insights compared to existing work \cite{grussler2021internally}; (iii) our earlier results on strict variation bounding have been published without proofs. These proofs can only be found in the present manuscript. 

The remainder of the paper is organized as follows. We begin with some extensive preliminaries in Section \ref{sec:prelim} that will enable us to present our main results in Section \ref{sec:main}. {In Section~\ref{sec:ex_pos_opt_ctr}, we relate our results to problem of bounding the number of over- and undershoots in a step response, establishing families of externally positive systems as well as bounding the number of bangs in bounded optimal control problems.} In Section~\ref{sec:examples}, our results are illustrated by examples and conclusions are drawn in Section \ref{sec:conclusion}. Proofs are stated in the Appendix.

\section{Preliminaries}
\label{sec:prelim}
In this section, we introduce basic notations and provide an extensive review of concepts that are essential for the presentation of our new results.

\subsection{Notations}
We write $\mathds{Z}$ for the set of integers and $\mathds{R}$ for the set of reals, with  $\mathds{Z}_{\ge 0}$ and $\mathds{R}_{\ge 0}$ standing for the respective subsets of nonnegative elements -- the corresponding notations are also used for subsets starting from non-zero values, strict inequality as well as reversed inequality signs. The set of real sequences with indices in $\mathbb{Z}$ is denoted by 
$\mathbb{R}^{\mathbb{Z}}$. For matrices $X = (x_{ij}) \in \Rmn$, we say that $X$ is
\emph{nonnegative},
$X \geq 0$ 
or $X \in \Rmn_{\geq 0}$
if all elements $x_{ij} \in \mathbb{R}_{\geq 0}$ -- corresponding notations are used for matrices with strictly positive entries and reversed inequality signs. These notations are also used for sequences $x = (x_i) \in \mathbb{R}^{\mathbb{Z}}$. For $k, l \in \mathds{Z}$, we write $(k:l) := \{k,k+1,\dots,l\}$, $k \leq l$. In the case $k > l$, the notation represents the empty set. If $X\in \Rnn$, then $\sigma(X) = \{\lambda_1(X),\dots,\lambda_n(X)\}$ denotes its \emph{spectrum}, where the eigenvalues are ordered by descending absolute value, i.e., $\lambda_1(X)$ is the eigenvalue with the largest magnitude, counting multiplicity. If the magnitude of two eigenvalues coincides, we sub-sort them by decreasing real part. The identity matrix in $\Rnn$ is denoted by $I_n$ or if the dimensions are obvious, simply by $I$. A \emph{(consecutive) $j$-minor} of $X \in \Rmn$ is a minor which is constructed of $j$ columns and $j$ rows of $X$ (with consecutive indices). The set of $j$-consecutive minors with column (row) indices $\nset{j}$ are said to be the $j$-\emph{column}(\emph{row}) \emph{initial minors}. The submatrix with rows $I \subset
(1:n)$ and columns $J \subset (1:m)$ is written as $\submatrix{X}{I}{J}$. In case of subvectors, we simply write $x_I$. With slight abuse of notation, we also use this to denote subsets of ordered index sets $x \in \mathcal{I}_{n,r}$, where 
\begin{equation*}
\mathcal{I}_{n,r} := \{ v = \{v_1,\dots,v_r\} \subset \mathds{N}: 1\leq v_1 < \dots < v_r \leq n \},
\end{equation*} 
and the \emph{complement} of $x$ is defined as the unique $\compset{x} \in \mathcal{I}_{n,n-r}$ such that $ \compset{x} \cup x = (1:n)$. Finally, if $X \in \Rnn$ is such that $\det \submatrix{X}{(2:n-1)}{(2:n-1)} \neq 0$, then the following identity exists \cite[Equation 0.8.11]{horn2012matrix}: 
\begin{multline}\label{eq:Desnanot_Jacobi_identity}
    \det \submatrix{X}{\nset{n}}{\nset{n}} \det \submatrix{X}{(2:n-1)}{(2:n-1)}\\
    =
    \det \submatrix{X}{\nset{n-1}}{\nset{n-1}} \det \submatrix{X}{(2:n)}{(2:n)}\\
    -
    \det \submatrix{X}{\nset{n-1}}{(2:n)} \det \submatrix{X}{(2:n)}{\nset{n-1}}
\end{multline}
\subsection{Variation diminishing maps}
\label{sec:vardim}
The \emph{variation} of a sequence or vector $u$ is defined 
as the number of sign changes in $u$. We employ two versions that only differ in the treatment of zero entries. 
\begin{equation*}
\vari{u} := \sum_{i} \mathbb{1}_{\mathbb{R}_{< 0}}(\tilde{u}_i \tilde{u}_{i+1}), 
    \quad 
    \vari{0} := -1
\end{equation*}
where $\tilde{u}$ is the vector resulting from deleting all zeros in $u$ and $\mathbb{1}_{\mathbb{A}}(x)$ is the indicator function with subset $\mathbb{A}$, i.e., $\mathbb{1}_{\mathbb{A}}(x) = 1$ if $x \in \mathbb{A}$ and zero otherwise. Further, {the \emph{strict variation} is defined by}
\begin{equation*}
\variz{u} := \sum_{i} \mathbb{1}_{\mathbb{R}_{< 0}}(\bar{u}_i \bar{u}_{i+1}), 
\end{equation*}
where $\bar{u}$ is the vector resulting from replacing zeros by elements that maximize the resulting sum. Obviously, $\vari{u} \leq \variz{u}$, {but equality does not necessarily need to hold, e.g., $\vari{(1 \; 0 \; 2)} = \vari{(1 \; 2)} = 0$, but $\variz{(1 \; 0 \; 2)} = \vari{(1 \; -1 \; 2)}  = 2$. }

Most essential to this work is the definition of variation bounding. 
\begin{defn}
\label{def:svb_k} \label{def:vb_k}
A linear map $u \mapsto X u$ is said to be \emph{$k$-variation bounding} ($\vb{k}$), 
$k \in \mathds{Z}_{\ge0}$, if for all $u\neq 0$ with $\vari{u} \leq k$ it holds that
\begin{equation*}
    \vari{Xu} \leq k 
\end{equation*}
If $\vari{Xu}$ can be replaced with $\variz{Xu}$, the mapping is said to be \emph{strictly $k$-variation bounding} ($\svb{k}$). For brevity, we simply say $X$ is $\text{(S)}\vb{k}$.
\end{defn}
A special case of variation bounding is that of (order-preserving) variation diminishment. 
\begin{defn}\label{def:ovd_k}
A linear operator $X$ is said to be \emph{$k$-variation diminishing} ($\vd{k}$), if $X$ is $\vb{j}$ for all $j \in \{0,\dots,k\}$. If additionally, the sign of the first non-zero {element} in $u$ and $Xu$, {respectively}, coincide whenever $\vari{u} = \vari{Xu}$, we say that $X$ is \emph{order-preserving $k$-variation diminishing} ($\ovd{k}$).
\end{defn}

\subsection{Total Positivity Theory}
\emph{Total positivity theory} \cite{karlin1968total} is known to provide algebraic conditions for the $\ovd{k-1}$ and $\svb{k-1}$ property by means of compound matrices. For $X \in \Rmn$, the $(i,j)$-th entry of the so-called \emph{r-th multiplicative compound matrix} 
$\compound{X}{r} \in \mathbb{R}^{\binom{n}{r} \times \binom{m}{r}}$ is defined via $\det(\submatrix{X}{I}{J})$, where $I$ and $J$ are the $i$-th and $j$-th element of the $r$-tuples in $\mathcal{I}_{n,r}$ and 
$\mathcal{I}_{m,r}$, respectively, assuming \emph{lexicographical ordering}.
For example, if $X \in \mathbb{R}^{3 \times 3}$, then
\begin{align*}
\begin{pmatrix}
\det(X_{\{1,2 \},\{1,2 \}}) & \det(X_{\{1,2 \},\{1,3\}}) & \det(X_{\{1,2 \},\{2,3\}})\\
\det(X_{\{1,3 \},\{1,2 \}}) & \det(X_{\{1,3 \},\{1,3\}}) & \det(X_{\{1,3 \},\{2,3\}})\\
\det(X_{\{2,3 \},\{1,2 \}}) & \det(X_{\{2,3 \},\{1,3\}}) & \det(X_{\{2,3 \},\{2,3\}})\\
\end{pmatrix}
\end{align*}
equals $\compound{X}{2}$. In our derivations, the following properties of the multiplicative compound matrix will be elementary (see, e.g., \cite[Section~6]{fiedler2008special} and \cite[Subsection~0.8.1]{horn2012matrix}).
\begin{lem}\label{lem:compound_mat}
	Let $X \in \mathbb{R}^{n \times p}$ and $Y \in \mathbb{R}^{p \times m}$.
	\begin{enumerate}[i)]
		\item $\compound{(XY)}{r} = \compound{X}{r}\compound{Y}{r}$ (Cauchy-Binet formula). \label{item:Cauchy_Binet}
		\item  For $ p = n $: $\sigma(\compound{X}{r}) = \{\prod_{i \in I} \lambda_i(X): I \in \mathcal{I}_{n,r} \}$. \label{item:compound_eigen}
  \item For $p=n$ and $\det(X) \neq 0$: $\compound{(X^{-1})}{r} = \compound{X}{r}^{-1}$ \label{item:compound_inv}.
  \item If $k = \rk(X)$, then $\rk(\compound{X}{k}) = 1$ \label{item:rk_k_compound_k}
	\end{enumerate} 
\end{lem}
It is readily seen that $X$ is $\ovd{0}$ if only if $X = \compound{X}{1} \geq 0$. This equivalence can be generalized to higher orders by the following definitions and characterizations (see~\cite[Prop.~7]{grussler2020variation} and \cite[Theorem 5.1.1]{karlin1968total}). 
\begin{defn}\label{def:k_pos_matrix}
    Let $X \in \Rmn$ and $k \leq \min\{m,n\}$. $X$ is called 
    \begin{enumerate}[i.]
        \item \emph{(strictly) $k$-sign consistent} ($\text{(S)}\sgc{k}$) if $\compound{X}{k} \geq (>) 0$ or $\compound{X}{k} \leq (<) 0$.
        \item \emph{(strictly) $k$-sign regular} ($\text{(S)}\sr{k}$) if $X$ is $\text{(S)}\sgc{j}$ for all $j \in (1:k)$. 
        \item \emph{(strictly) $k$-positive} if $\compound{X}{j} \geq (>) 0$ for all $j \in (1:k)$. In case of $k = \min\{m,n\}$, $X$ is also called \emph{(strictly) totally 
    positive}.
        \end{enumerate}

\end{defn}
\begin{prop} \label{prop:k_pos_mat_var}
	Let $X \in \Rmn$ with $n \geq m$. Then, $X$ is $\tp{k}$ with $k \in (1:m)$ if and only if $X$ is $\ovd{k-1}$.
\end{prop}
\begin{prop} \label{prop:n_ssc_mat_svb}
	$X \in \Rmn$, $n > m$, is $\ssgc{m}$ if and only if $X$ is $\svb{m-1}$.
\end{prop}
Note that one is only interested in cases with $n > m$, because any $X \in \mathds{R}^{n \times n}$ is $\svb{n-1}$. In order to extend strict to the non-strict cases, the following proposition (see, e.g., \cite[Propositions~5.1.1 \& 5.1.2]{karlin1968total}, \cite[Propositions~2.7]{grussler2021internally}) and lemma \cite[Lemma~5.1.1]{karlin1968total} will be used.  %
\begin{prop} \label{prop:T_sigma}
Let $T(\sigma) \in \mathds{R}^{n\times n}$ be given by $\submatrix{T(\sigma)}{i}{j}=e^{-\sigma (i-j)^2}$, with $\sigma >0$, and let $X \in \Rmn$ with $m\leq n$. Then for $k \leq m$, the following hold:
\begin{enumerate}[i.]
    \item $T(\sigma)$ is strictly totally positive 
    \item $T(\sigma) \rightarrow I$ as $\sigma \rightarrow \infty$, and $T(\sigma)X \rightarrow X$ as $\sigma \rightarrow \infty$
    \item if $X$ is $\sgc{k}$, and if all $k$ columns of $X$ are linearly independent, then $T(\sigma) X$ is $\ssgc{k}$. \label{item:T_sigma_X_SSC_k}
    \item if $T(\sigma) X$ is $\sgc{k}$ for all $\sigma > 0$, then $X$ is $\sgc{k}$. \label{item:X_SC_k}
    \item if $X$ is $\vb{k-1}$ and $\rk(X) = m$, then $T(\sigma) X$ is $\svb{k-1}$ \label{item:T_sigma_X_SVB_k}
\end{enumerate}
\end{prop}
\begin{lem} \label{lem:var_lim}
Let $x: \mathds{R} \to \mathds{R}^n$ be such that $x(\sigma) \rightarrow x$ as $\sigma \rightarrow \infty$. Then, $\vari{x} \leq \lim_{\sigma \rightarrow \infty} \variz{x(\sigma)}$ and $\lim_{\sigma \rightarrow \infty} \vari{x(\sigma)} \leq \variz{x}$.
\end{lem}
For example, Proposition~\ref{prop:T_sigma} in conjunction with Lemma~\ref{lem:var_lim} characterizes $\vb{m-1}$ by $\sgc{m}$ as in \cite[Theorems 5.1.3 \& 5.1.3']{karlin1968total}.
\begin{prop}\label{prop:sc_k_mat_vb_m}
	Let $X \in \Rmn$, $n > m$ and $\mathrm{rank}\; m$. $X$ is $\sgc{m}$ if and only if $X$ is $\vb{m-1}$.
\end{prop}
In this work, we will extend {this} characterization to the case of $\text{(S)}\vb{k}$. 

Next, we will review how checking the sign of the elements in $\compound{X}{k}$ can be simplified compared to computing all of its entries. This does will not only make $\text{(S)}\vb{k-1}$ more tractable, but is particularly important for our applications to infinite-dimensional operators. 
We begin with a sufficient condition for $m$-sign consistency \cite[Theorem 2.3.1]{karlin1968total}, which only requires to consider consecutive row and column initial minors.
\begin{prop}\label{prop:col_intial_m_strict_sign_consist}
    Let $X \in \Rmn, \; n\geq m$ be such that
    \begin{enumerate}[i.]
        \item the $k$-column initial minors of $X$ are positive $k = 1,2,..m-1$
        \item the $m$ column initial minors are nonnegative (positive).
    \end{enumerate}
    Then all $m$-order minors of $X$ are nonnegative (positive). The same is true for $m \geq n$ by replacing column with row initial minors.
\end{prop}
An application of this result
yields the following (sufficient) certificate (see~\cite[Proposition~8]{grussler2020variation}) for (non-)strict $k$-positivity.  
\begin{prop}
    \label{prop:consecutive_old}
	Let $X \in \Rmn$, $k \leq \min \{n,m\}$, be such that
	\begin{enumerate}[i.]
	    \item all consecutive $r$-minors of $X$ are 
	    positive, $r \in (1:k-1)$,
	    \item all consecutive $k$-minors of $X$ are
	    nonnegative (positive).
	\end{enumerate}
	Then, $X$ is (strictly) \tp{k}. In the strict case, these are also necessary conditions.  
\end{prop}
 A related result is the following sufficient and necessary condition for $\ssgc{m}$ \cite[Theorem 2.2]{pena_matrices_1995} and its application to $\ssgc{k}$ case.
\begin{prop}\label{prop:k_ssc__mat_pena_i}
    $X \in \Rmn$, $n>m$, is $\ssgc{m}$ if and only if the minors $\det(\submatrix{X}{\alpha}{\nset{m}})$ have the same strict sign for all $\alpha = \{\nset{m-r},(t:t+r-1)\}$ with $1\leq r \leq m $ and $m-r \leq t \leq n-r+1$.
\end{prop}
For example, checking whether
\begin{equation*}
    X = \begin{pmatrix}
        1 & 1 \\ 
        1 & 2 \\ 
        1 & 3 \\
        1 & 4
    \end{pmatrix} \in \mathds{R}^{4 \times 2}
\end{equation*}
is $\ssgc{2}$ only requires to verify that
\begin{align*}
    r=1:& \quad \det{\begin{pmatrix} 1 & 1\\ 1 & 2 \end{pmatrix}},\;
\det{\begin{pmatrix} 1 & 1 \\ 1 & 3 \end{pmatrix}}, \;
\det{\begin{pmatrix} 1 & 1 \\ 1 & 4 \end{pmatrix}}
    \\
    r=2:& \quad 
    \det{\begin{pmatrix} 1 & 1\\ 1 & 2 \end{pmatrix}}, \;
    \det{\begin{pmatrix} 1 & 2 \\ 1 & 3 \end{pmatrix}}, \;
    \det{\begin{pmatrix} 1 & 3\\ 1 & 4 \end{pmatrix}}.
\end{align*}
have the same strict sign. Propositions~\ref{prop:col_intial_m_strict_sign_consist} -- \ref{prop:k_ssc__mat_pena_i} turn the combinatorial complexity of verifying $\ssgc{m}$ or $\text{(S)TP}_k$ via $\compound{X}{m}$ into polynomial complexity.

A simple application of Propositions~ \ref{prop:k_ssc__mat_pena_i} to all $k$ columns of $X$ also provides a characterization of $\ssgc{k}$.
\begin{cor}\label{cor:k_ssc__mat_pena_ii}
     $X \in \Rmn$, $n>m$, is $\ssgc{k}$, $k \leq m$ if and only if the minors $\det (\submatrix{X}{\alpha}{\beta})$ have the same strict sign for all $\alpha = \{\nset{k-r},(t:t+r-1)\}$ with $r \in (1:k)$ and $t \in (k-r+1 : n-r+1)$ and all $\beta = \{\nset{k-\bar r},(\bar t:\bar t+ \bar r-1)\}$ with $\bar r \in (1:k)$ and $\bar t \in (k- \bar r + 1 : m- \bar r+1)$.
\end{cor}
Proposition~\ref{prop:k_ssc__mat_pena_i} utilizes a bijection between the $m$-order minors of a matrix $X \in \Rmn$ and the minors of a related transformed matrix \cite{pena_matrices_1995}, which we will also use for a new non-strict counter-parts. 
\begin{lem}\label{lem:pena_bijection}
    Let $X \in \Rmn, \; n>m$, such that $\det \submatrix{X}{\nset{m}}{\nset{m}} \ne 0$ and 
    \begin{equation*}
       C :=  \submatrix{X}{(n-m+1:n)}{\nset{m}} \left(\submatrix{X}{\nset{m}}{\nset{m}}\right)^{-1} K, 
    \end{equation*}
    where $K \in \mathds{R}^{m \times m}$ is defined by
    \begin{equation*}
        k_{ij} = 
        \begin{cases}
            (-1)^{j-1} & i+j = m+1\\
            0 & \text{otherwise}
        \end{cases}.
    \end{equation*}    
    Then there is a bijection between the $m$-order minors $\det \submatrix{X}{\gamma}{\nset{m}}$ with $\gamma \in \mathcal{I}_{n,m} \backslash \nset{m}$ and all minors of $C$. Concretely, $\det \submatrix{C}{\alpha}{\beta}, \; \alpha \in \mathcal{I}_{n-m,r}, \; \beta \in \mathcal{I}_{m,r}$ for all $r \in \nset{\min \{m,n-m\}}$, where
    \begin{align*}
        &\gamma_{m-r+1-j} = m + 1 - \compset{\beta}_j &\text{for } j\in \nset{m-r}\\
        &\gamma_{m-r+i} = m + \alpha_i & \text{for } i\in \nset{r}
    \end{align*}
    and $r \in \nset{\min\{m,n-m\}}$ is the order of the respective minor of $C$. The sign of the minors are equal up to the sign of $\det X_{(1:m),(1:m)}$.
\end{lem}

\subsection{Variation diminishing observability operators}
The \emph{observability operator} of \eqref{eq:SISO_d}, {$\Obs{}(A,c): \mathds{R}^n \to \mathds{R}^{\mathds{Z}_{\geq 0}}$} is defined by \begin{align}
 	(\Obs{}(A,c) x_0)(t) &:= cA^t x_0, \ x_0 \in \mathds{R}^n,  \ t \geq 0. \label{eq:obsv_op}
 \end{align}
 As discussed in \cite[Lemma~3.4]{grussler2021internally}, the $\ovd{k}$ property of 
$\Obs{}(A,c)$ is equivalent to
\begin{equation*}
    \Obs{t}(A,c)^\transp := \begin{pmatrix}
        c^\transp & A^\transp c^\transp & \dots & {A^{t-1}}^\transp c^\transp 
    \end{pmatrix}
\end{equation*}
being $\ovd{k}$ for all $t\geq k$. For brevity, we will mostly drop $(A,c)$ and write $\Obs{}$ and $\Obs{t}$, respectively. A sufficient condition based on the assumption of a $k$-positive $A$ was shown in \cite[Theorem~3.5]{grussler2021internally}. 
\begin{prop} \label{prop:k_pos_obsv}
 If $A$ is $k$-positive and $\compound{\Obs{j}}{j} \geq 0$ for $1\leq j \leq k$, then $\Obs{t}$ is $k$-positive for all $t \geq k$.
\end{prop}
Finally, we will abbreviate a discrete-time LTI system (\ref{eq:SISO_d}) by the triple $(A,b,c)$ and say that it is \emph{strictly externally positive} (\ref{eq:SISO_d}) if $g(t) = cA^{t-1}b > 0$ for all $t \geq 1$. Correspondingly, we use the term \emph{strictly externally negative} if the inequality is reversed.

\section{Main Results}
\label{sec:main}
The main goal of this work is to verify the general case of $\Obs{}$ being $(S)\vb{k-1}$ using $(S)\sgc{k}$. This enables us to derive numerical certificates as well as necessary conditions in terms of the eigenvalues of $A$. In particular, we also get necessary and sufficient conditions for $(O)\vd{k}$ $\Obs{}$ and as such no longer require $k$-positivity of $A$ as in Proposition \ref{prop:k_pos_obsv}. Finally, by substitution of $(A,c)$ with $(A^\transp,b^\transp)$, these results may also be used to checked if the \emph{controllability operator} $\mathcal{C}(A,b) := \Obs{}(A^\transp,b^\transp)^\transp$ of (\ref{eq:SISO_d}) is $(S)\vb{k}$ and as result provide a first sufficient certificate for the Hankel operator $\mathcal{H}_g = \mathcal{O}(A,c)\mathcal{C}(A,b)$ of (\ref{eq:SISO_d}) to be $(S)\vb{k}$. 

Our investigations start with the $\svb{k-1}$ case and are subsequently extended to the $\vb{k-1}$ case. 

\subsection{The $\svb{k-1}$ case}
 Analogously to the $\ovd{k-1}$ case, one can show that $\Obs{}$ is $\svb{k-1}$ if and only if $\Obs{t}$ is $\svb{k-1}$ for all $t \geq k-1$. Thus, by Proposition \ref{prop:n_ssc_mat_svb}, it suffices to check $\ssgc{n-1}$ of $\Obs{t}$ for all $t\geq n-1$. Using Proposition \ref{prop:k_ssc__mat_pena_i}, this is equivalent to the strict positivity/negativity of certain sequences of $n-1$-minors. Next, we will show that these sequences of $n-1$- minors correspond to impulse responses of related LTI system. Checking $\svb{k-1}$ then becomes equivalent to the verification of strict external positivity of $n$ LTI systems for which numerically efficient certificates can be found, e.g., in \cite{grussler2019tractable,taghavian2023external}. 

\begin{thm} 
	\label{thm:main1}
	Let $(A,c)$ be observable. Then, $\Obs{}$ is $\svb{n-1}$ if and only if $(\tilde{A}_r,\tilde{b}_r,\tilde{c}_r)$ is strictly externally positive for all $r \in (1:n)$ with $\tilde A_r = \compound{A}{r}$, $\tilde c_r = \pcompound{\Obs{r}}{r}$ and 
	\begin{align*}
		\tilde b_r = \pcompound{A^{n-r}\left(\Obs{n}\right)^{-1} \begin{pmatrix} 0 \\ I_r  \end{pmatrix}}{r}.
	\end{align*}
\end{thm}
In order to proceed with the general $\svb{k-1}$ case, we need the following generalization of Proposition \ref{prop:n_ssc_mat_svb}, which has appeared earlier in several related forms (see, e.g., \cite{karlin1968total,weiss2019generalization}), but to the best of our knowledge not been stated and derived in given full generality.
\begin{prop}\label{prop:ssc_k_mat_svb_k}
	Let $X \in \Rmn$ and $k < \min(m+1,n)$. Then, $X$ is $\ssgc{k}$ if and only if $X$ is $\svb{k-1}$.
\end{prop}
 With the next theorem we achieve a method for checking whether $\Obs{}$ is $\svb{k-1}$.  
\begin{thm}\label{thm:main2}
Let $(A,c)$ be observable, $k \in (1:n)$, $r \in (1:k)$ and $\beta \in \mathcal{I}_{n,k}$. Then, for 
$N\in \mathds{N}_{\geq 1}, \; t \in \nset{N-k+1}$ and $\alpha = \{\nset{k-r},(k-r+t:k+t-1)\}$, it holds that 
\begin{equation*}
    \det(\mathcal{O}^{N}_{\alpha,\beta}) = \tilde{c}_r \tilde{A}_r^{t-1}\tilde{b}_{k,r,\beta},
\end{equation*}
where $\tilde A_{r} := \compound{A}{r}$, $\tilde{c}_r := \compound{\Obs{r}}{r}$ and 
	\begin{align*}
		\tilde b_{k,r,\beta} &:= 
		\begin{pmatrix} \pcompound{A^{k-r}\left(\Obs{n}\right)^{-1}\begin{pmatrix} 0 \\ I_{n-k+r} \end{pmatrix}}{r} & 0 \end{pmatrix} \pcompound{\Obs{n}P_{\beta}}{k}
	\end{align*}
with $P_\beta := \submatrix{I}{\nset{n}}{\beta}$. Therefore, $\Obs{}$ is $\svb{k}$ if and only if $(\tilde{A}_{r},\tilde{b}_{k,r,\beta},\tilde{c}_r)$ is strictly externally positive/negative for all $(r,\beta)$ with $r \in (1:k)$ and $\beta$ as in Corollary \ref{cor:k_ssc__mat_pena_ii}.%
\end{thm}
Since strict externally positivity requires a dominant positive pole (see, e.g., \cite{ohta1984reachability,farina2011positive,benvenuti2004tutorial}), the following eigenvalue constraint can be shown as a consequence of Item~\ref{item:compound_eigen} in Lemma~\ref{lem:compound_mat}.
\begin{cor} \label{cor:eigen}
{Let $(A,c)$ be observable such that $\mathcal{O}$ is $\ssgc{k}$. Then,  $\lambda_1(A),\dots,\lambda_k(A)$ are simple and strictly positive.}
\end{cor}
Interestingly, this necessary condition is also implied by $k$-positivity of $A$ as in Proposition~\ref{prop:k_pos_obsv} as well as by the $k$-positivity of the Hankel and Toeplitz operators \cite{grussler2020variation}. In particular, if $k = n$, then there exists a $T \in \mathds{R}^{n \times n}$ such that $T^{-1}AT = \diag(\lambda_1(A),\dots,\lambda_n(A)) \geq 0$ and $cT = \begin{pmatrix}
    1 & \cdots & 1
\end{pmatrix}$, which as shown in \cite{grussler2021internally} fulfills the requirements of Proposition~\ref{prop:k_pos_obsv}.

 \subsection{The $\vb{k-1}$ case}
Next, we will extend the above strict version to their corresponding non-strict counterparts. To this end, we will start by relaxing Proposition \ref{prop:ssc_k_mat_svb_k} and Corollary \ref{cor:k_ssc__mat_pena_ii}.
Unlike the strict case, we need to divide the characterization of $\vb{k-1}$ matrices into two distinct cases.
\begin{thm}\label{thm:sc_k_mat_vb_k_rank_k}
	Let $X \in \Rmn$, $n > m$ and $k= \rk\; X <m$. Then the following are equivalent
    \begin{enumerate}[i.]
        \item Each column of $X_{[k]}$ has elements of the same sign. \label{item:column_sgc}
        \item For every $u\in \mathds{R}^m$, $\vari{Xu} \leq k-1$.  \label{item:sign_bound}
    \end{enumerate}
\end{thm}
\begin{thm} \label{thm:sc_k_mat_vb_k_rank_gt_k}
    Let $X \in \Rmn$ and $k < \rk (X)$ be such that any $k$ columns of $X$ are linearly independent. Then 
    \begin{equation*}
        X \; \text{is} \; \vb{k-1} \iff X \; \text{is} \; \sgc{k}
    \end{equation*}
\end{thm}
Proposition~\ref{prop:k_ssc__mat_pena_i} is a result of 
Lemma \ref{lem:pena_bijection} in conjunction with the following characterization of strictly totally positive matrices (see, e.g., \cite[Theorem 2.1]{fallat2017total}). 
\begin{prop}\label{prop:row_col_initial_minors_strict_tot_pos}
     $X \in \Rmn$  is strictly totally positive if and only if all row and column initial minors of $X$ are positive.
\end{prop}
For completion, a new proof of Proposition~\ref{prop:row_col_initial_minors_strict_tot_pos}, similar to \cite[Theorem~2.3]{pinkus2009totally}, is stated in the appendix. Based on that, we show the following new non-strict relaxation similar to Proposition~\ref{prop:consecutive_old}.
\begin{prop}\label{prop:row_col_initial_minors_tot_pos}
     Let $X \in \Rmn$ be such that row and column initial $k$-order minors of $X$ are positive for $k = 1, 2, ..., \min\{n-1,m-1\}$ and nonnegative for $k = \min\{n,m\}$. Then, $X$ is totally positive and all minors of order less than $\min\{n,m\}$ are positive. 
\end{prop}
Non-strict analogues of Proposition \ref{prop:k_ssc__mat_pena_i} and Corollary~\ref{cor:k_ssc__mat_pena_ii} then read as follows. 
\begin{thm}\label{thm:k_sc__mat_pena_i}
    $X \in \Rmn$, $n \geq 2m$, is $\sgc{m}$ if the minors $\det(\submatrix{X}{\alpha}{\nset{m}})$ have the same strict sign for all $\alpha = \{\nset{m-r}, (t:t+r-1)\}$ with $1\leq r < m,\;\;m-r \leq t \leq n-r+1$, and $\alpha = (t:t+m-1),\;\; 1 \leq t \leq m$. In addition the minors need to have the same non-strict sign for $\alpha = (t:t+m-1),\;\; m+1 \leq t \leq n-m+1$.
\end{thm}

\begin{cor}\label{cor:k_sc__mat_pena_ii}
     $X \in \Rmn$, $n > m$, is $\sgc{k}$, $ 2 k \leq m$ if the minors $\det (\submatrix{X}{\alpha}{\beta})$ have the same strict sign for all $\alpha = \{\nset{k-r},(t:t+r-1)\}$ with $r \in (1:k)$ and $t \in (k-r+1 : n-r+1)$ and all $\beta = \{\nset{k-\bar r},(\bar t:\bar t+ \bar r-1)\}$ with $\bar r \in (1:k)$ and $\bar t \in (k- \bar r + 1 : m- \bar r+1)$, except for the cases where $\alpha = (t:t+r-1)$ $r=k$, $t \in (k+1 : n-k+1)$ and $\beta = (\bar t:\bar t+ \bar r-1)$, $\bar r = k$, $\bar t \in (k+1 : m- k +1)$, where a non-strict sign is sufficient.
\end{cor}
Equipped with our new non-strict matrix conditions, we are ready to state the following sufficient analogues of Theorems \ref{thm:main1} \&\ref{thm:main2}. 
\begin{prop} 
	\label{prop:nonstrict_main1}
	Let $(A,c)$ be observable and $(\tilde{A}_r,\tilde{b}_r,\tilde{c}_r)$ be defined by $\tilde A_r := \compound{A}{r}$, $\tilde c_r := \pcompound{\Obs{r}}{r}$ and 
	\begin{align*}
		\tilde b_r := \pcompound{A^{n-r}\left(\Obs{n}\right)^{-1} \begin{pmatrix} 0 \\ I_r  \end{pmatrix}}{r}.
	\end{align*}
 for all $r \in (1:n)$. Then, $\Obs{}$ is $\vb{n-1}$ {if} $(\tilde{A}_r,\tilde{b}_r,\tilde{c}_r)$ is strictly externally positive for all $r \in (1:n-1)$ and externally positive for $r=n$ with $\tilde{c}_n\tilde{A}^{l-1}_n\tilde{b} > 0$ for all $l \in (1:n-1)$. 
\end{prop}

\begin{prop}\label{prop:nonstrict_main2}
Let $(A,c)$ be observable, $k \in (1:n)$, $r \in (1:k)$ and $\beta \in \mathcal{I}_{n,k}$. Then, for 
$N\in \mathds{N}_{\geq 1}, \; t \in \nset{N-k+1}$ and $\alpha = \{\nset{k-r},(k-r+t:k+t-1)\}$, it holds that 
\begin{equation*}
    \det(\mathcal{O}^{N}_{\alpha,\beta}) = \tilde{c}_r \tilde{A}_r^{t-1}\tilde{b}_{k,r,\beta},
\end{equation*}
where $\tilde A_{r} := \compound{A}{r}$, $\tilde{c}_r := \compound{\Obs{r}}{r}$ and 
	\begin{align*}
		\tilde b_{k,r,\beta} &:= 
		\begin{pmatrix} \pcompound{A^{k-r}\left(\Obs{n}\right)^{-1}\begin{pmatrix} 0 \\ I_{n-k+r} \end{pmatrix}}{r} & 0 \end{pmatrix} \pcompound{\Obs{n}P_{\beta}}{k}
	\end{align*}
with $P_\beta := \submatrix{I}{\nset{n}}{\beta}$. Therefore, $\Obs{}$ is $\vb{k-1}$ if  $(\tilde{A}_{r},\tilde{b}_{k,r,\beta},\tilde{c}_r)$ is strictly externally positive/negative for all $(r,\beta)$ with $r \in (1:k)$ and $\beta$ as in Corollary \ref{cor:k_ssc__mat_pena_ii}, except for the following cases $\beta = (\bar t:\bar t+ \bar r-1)$, $\bar r = k$, $\bar t \in (k+1 : m- k +1)$, where non-strict external positivity/negativity is sufficient if $\tilde{c}_r\tilde{A}_{[r]}^{l-1}\tilde b_{k,r,\beta} \neq 0$ for all $l \in (1:k)$.
\end{prop}
 {Note that in contrast to Theorems ~\ref{thm:main1} and \ref{thm:main2}, these results are only sufficient, because Theorem~\ref{thm:sc_k_mat_vb_k_rank_gt_k} and Corollary~\ref{cor:k_sc__mat_pena_ii} are generally not necessary for matrices to be $\sgc{k}$.}
\subsection{The $\ovd{k}\slash \vd{k}$ case}
An application of Theorem~\ref{thm:main2} also provides a new method to check $\ovd{k}$. To this end note that the systems in Theorem~\ref{thm:main2}, which represent the $k$-consecutive minors of $\Obs{}$ are the pairs with $r = k$. In conjunction with Proposition \ref{prop:consecutive_old}, this provides the following corollary for $k$-positive $\Obs{}$.
\begin{cor}\label{cor:k_pos}
	Let $(A,c)$ be observable. Using the notation of Theorem \ref{thm:main2}, it holds that if
	\begin{enumerate}[i.]
		\item $(\tilde{A}_r,b_{r,r,\beta},\tilde{c}_r)$ is strictly externally positive, $r \in (1:k-1)$
		\item {$(\tilde{A}_k,b_{k,k,\beta},\tilde{c}_k)$} is (strictly) externally positive,
	\end{enumerate}
	then $\Obs{}$ is (strictly) $k$-positive.
\end{cor}
It is important to note that unlike Proposition \ref{prop:k_pos_obsv}, Corollary \ref{cor:k_pos} does not require $A$ to be $k$-positive. In Subsection~\ref{subsec:A_not_kpos}, we even present an example where the system does not permit any realization with $k$-positive $A$.

Finally, for $\vd{k}$ there exists an equivalence to $\sr{k}$ with a full rank requirement \cite[Theorem 5.1.4]{karlin1968total}. Based on Theorem \ref{thm:sc_k_mat_vb_k_rank_gt_k}, we are able to relax this requirement, too. 

\begin{cor} \label{cor:sr_k_mat_vd_k_rank_gt_k}
    Let $X \in \Rmn$ and $\rk(X) > k$ be such that any $k$ columns of $X$ are linearly independent. Then 
    \begin{equation*}
        X \; \text{is} \; \vd{k-1} \iff X \; \text{is} \; \sr{k}
    \end{equation*}
\end{cor}

\section{Discussion}
\label{sec:ex_pos_opt_ctr}
In the following, it is discussed how our results relate to various problems in systems \& control. We start by discussing the case of bounding the number of sign changes of an impulse response of (\ref{eq:SISO_d}): 
\begin{enumerate}[I.]
	\item The number of under- and overshoots in a step response, which may be used as measure for the trade-off between fast rise vs. fast settling time, equals the number of sign changes of the impulse response. Yet, mathematically, little is known about the tractability of this measure. While there exist several lower bounds \cite{damm2014zero,swaroop1996some,elkhoury1993discrete}, aside from external positivity certificates \cite{altafini2016minimal,grussler2019tractable,vandeneijnden2025externally,taghavian2023external,farina2011positive}, only few methods provide generic upper bounds \cite{elkhoury1993discrete,grussler2021internally}. By noticing that  the impulse response of (\ref{eq:SISO_d}) can be expressed as $g = \Obs{}(A,c)b$, our results can be used to bound the number of sign changes in $g$ via
	\begin{equation}
		\vari{g} = \vari{\Obs{}(A,c)b} \leq k \label{eq:impulse_bound}
	\end{equation}
	whenever $\Obs{}(A,c)$ is $\vb{k}$ and $\vari{b} \leq k$. In particular, Corollary~\ref{cor:k_pos} can be used to check whether $\Obs{}(A,c)$ is $\ovd{1}$, which facilitates the following certificate for a family of externally positive systems:
	\begin{prop}
		Let $(A,b,c)$ be such that $\Obs{}(A,c)$ be $2$-positive, $\vari{b} \leq 1$. Further, assume that the first non-zero element in $b$ is negative and the first non-zero element of the impulse response of $(A,b,c)$ is positive. Then, $(A,b,c)$ is externally positive. 
		\end{prop}
		\vspace*{0.5 cm}
	\item Bounded optimal control problems such as 
\begin{equation}\label{eq:ctrl}
	\begin{aligned}
		& \underset{u}{\text{minimize}}
		& & \|u\|_{\ell_\infty} \\
		& \text{subject to}
		& & x({t+1}) = Ax(t)+bu(t), \\ 
		& & & x(0) = x_0,\; x(N) = 0
	\end{aligned}
\end{equation}
where $\| u \|_{\ell_\infty} := \max_{t \in (0:N-1)} |u(t)|$, are known to lead (e.g., by the alignment property \cite[Chapter~5]{luenberger1968optimization}) to solutions of the form
\begin{equation*}
	u^\ast(t) = u_{\max} \sign(\mu^\transp A^{-t-1}b) \text{ if } \mu^\transp A^{-t-1} \neq 0
\end{equation*}
for some $\mu \in \mathds{R}^{n}$, $u_{\max} > 0$. Otherwise, $u^\ast(t)$ may take any value in $[-u_{max},u_{\max}]$. In other words, the number of sign changes (often called \emph{bangs}) of $u^\ast$ are bounded by the strict variation of the impulse response $g_\mu(t) := \mu^\transp A^{-t} A^{-1} b$. Using Theorem~\ref{thm:main2}, we can conclude similar to \eqref{eq:impulse_bound} that if $\Con{}(A^{-1},A^{-1}b)$ is $\svb{k}$ and $\vari{\mu} \leq k$, then $u^\ast$ has at most $k$ bangs. This is an extension and discrete-time analog of the well-known continuous-time characterization, which asserts that if $(A,b)$ is controllable and $\sigma(A) \subset \mathds{R}$, then $u^\ast$ has at most $n-1$ bangs \cite{glad2000control}.
\end{enumerate}
Further, note that Theorems~\ref{thm:sc_k_mat_vb_k_rank_gt_k} and \ref{thm:k_sc__mat_pena_i} as well as  Propositions~\ref{prop:ssc_k_mat_svb_k} and \ref{prop:nonstrict_main1} have been instrumental for new tractable approaches in analyzing basis pursuit in structured problems, incl. fuel optimal control \cite{marmary2025tractable}, as well as the study of self-oscillations in relay feedback systems \cite{tong2025selfsustained,grussler2025discrete}.

\section{Illustrative Examples}
\label{sec:examples}
In this section, we want to illustrate our results based on three system examples \eqref{eq:SISO_d} for which we will bound the variation of their impulse responses via the relationship $g = \Obs{} b$. The three systems have the following distinct properties:
\begin{enumerate}[i.]
    \item $\Obs{}$ is $\ovd{1}$, but $A$ is not $2$-positive, i.e., Proposition \ref{prop:k_pos_obsv} does not apply for the given realization.
    \item $\Obs{}$ is $\svb{1}$, but not $2$-positive, i.e., variation diminishing results cannot be used, but our new results can be used to verify the variation bounding property. 
    \item $\Obs{}$ is $\vd{2}$, but there exists no other realization with $\sgc{2}$ $A$. As such, even a possible extension of Proposition~\ref{prop:k_pos_obsv} to $\sgc{2}$  would not be applicable. 
\end{enumerate}
{In all examples, we employ \cite{grussler2019tractable} to check external positivity/negativity of the impulse response as required by Theorem~\ref{thm:main2} and Corollary~\ref{cor:k_pos}.}
\subsection{Example 1: $\ovd{1}$ $\Obs{}$}
\begin{figure}
	\centering
	\begin{tikzpicture}
		\begin{axis}[height=4.5cm,
			width=\columnwidth,
			axis y line = left,
			axis x line = center,
			ylabel={$\tilde{g}_{1,1,\beta}(t)$}, 
			xlabel={$t$}, 
			xlabel style={right},
			xmax=10,
			]			
			\addplot+[ycomb,black,thick]
			table[x index=0,y index=1] 
			{g3_j1_r1_beta_1.txt};	\label{symbol:r1beta12_g3}
			\addplot+[ycomb,black,thick]
			table[x index=0,y index=1] 
			{g3_j1_r1_beta_2.txt};	\label{symbol:r1beta13_g3}
			\addplot+[ycomb,black,thick]
			table[x index=0,y index=1] 
			{g3_j1_r1_beta_3.txt};	\label{symbol:r1beta23_g3}
		\end{axis}
	\end{tikzpicture}
	\caption{Impulse responses $\tilde{g}_{1,1,\beta}(t)$ of $(\tilde{A}_1,b_{1,1,\beta},\tilde{c}_1)$ in Theorem \ref{thm:main2}:  \ref{symbol:r1beta12_g3} $\tilde{g}_{1,1,\{1\}}(t)$, \ref{symbol:r1beta13_g3} $\tilde{g}_{1,1,\{2\}}(t)$ and \ref{symbol:r1beta23_g3} $\tilde{g}_{1,1,\{3\}}(t)$ are strictly positive.}
	\label{fig:2pos_impulse_r1_j1}
\end{figure}
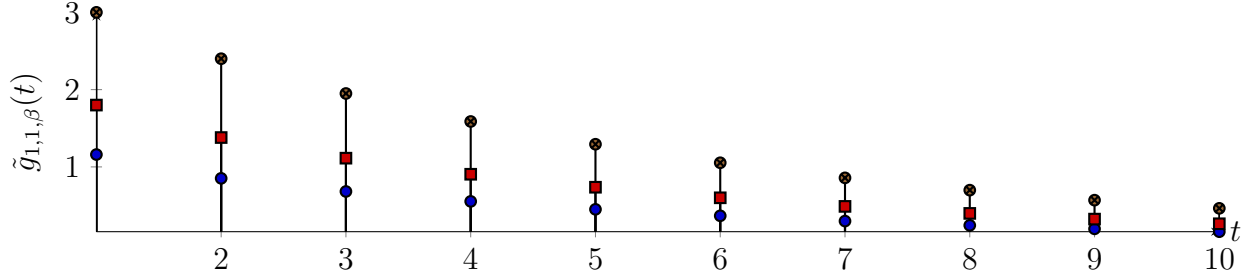

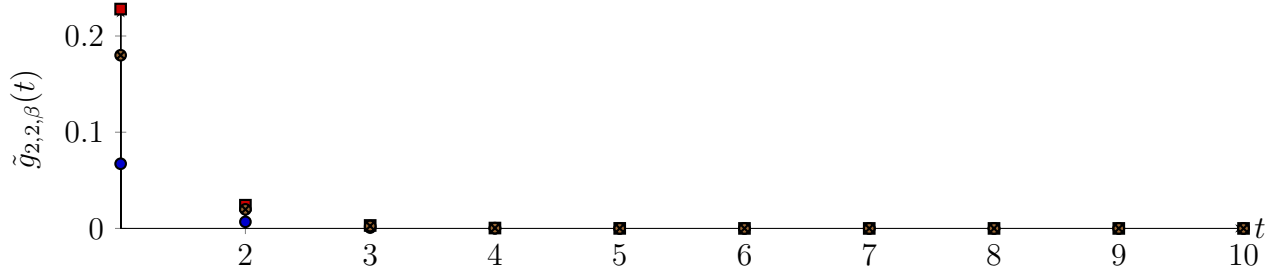
\begin{figure}
	\centering
	\begin{tikzpicture}
		\begin{axis}[height=4.5cm,
			width=\columnwidth,
			axis y line = left,
			axis x line = center,
			ylabel={$\tilde{g}_{2,2,\beta}(t)$}, 
			xlabel={$t$}, 
			xlabel style={right},
			xmax=10,
			]			
			\addplot+[ycomb,black,thick]
			table[x index=0,y index=1] 
			{g3_j2_r2_beta_12.txt};	\label{symbol:r2beta12_g3}
			\addplot+[ycomb,black,thick]
			table[x index=0,y index=1] 
			{g3_j2_r2_beta_13.txt};	\label{symbol:r2beta13_g3}
			\addplot+[ycomb,black,thick]
			table[x index=0,y index=1] 
			{g3_j2_r2_beta_23.txt};	\label{symbol:r2beta23_g3}
		\end{axis}
	\end{tikzpicture}
	\caption{Impulse responses $\tilde{g}_{2,2,\beta}(t)$ of $(\tilde{A}_2,b_{2,2,\beta},\tilde{c}_2)$ in Theorem \ref{thm:main2}:  \ref{symbol:r2beta12_g3} $\tilde{g}_{2,2,\{1,2\}}(t)$, \ref{symbol:r2beta13_g3} $\tilde{g}_{2,2,\{1,3\}}(t)$ and \ref{symbol:r2beta23_g3} $\tilde{g}_{2,2,\{2,3\}}(t)$ are strictly positive.}
	\label{fig:2pos_impulse_r2_j2}
\end{figure}

We begin by considering
\begin{align*}
	A = \begin{pmatrix}
		-1.20  & -1.50 &-1.88\\
		1.51  &  1.75  &  1.88\\
		-0.16  & -0.01 &   0.40
	\end{pmatrix}, \quad c^\transp = \begin{pmatrix}
		1.16  \\  1.8\\    3
	\end{pmatrix}.
\end{align*}
Since $A \not \geq 0$, $A$ is not $2$-positive and Proposition \ref{prop:k_pos_obsv} is not applicable. However, by checking external positivity of the impulse responses corresponding to consecutive minors of $\Obs{}$ (see~Figures \ref{fig:2pos_impulse_r1_j1} \& \ref{fig:2pos_impulse_r2_j2}), it follows from Corollary \ref{cor:k_pos} that $\Obs{}$ is $2$-positive. Thus, for all $b \in \mathds{R}^3$ with $\vari{b} \leq 1$, it follows that $\vari{g} \leq \variz{g} \leq 1$. 

\subsection{Example 2: $\svb{1}$ $\Obs{}$}
\begin{figure}
	\centering
	\begin{tikzpicture}
		\begin{axis}[height=4.5cm,
			width=\columnwidth,
			axis y line = left,
			axis x line = center,
			ylabel={$\tilde{g}_{2,1,\beta}(t)$}, 
			xlabel={$t$}, 
			xlabel style={right},
			xmax=10,
			]			
			\addplot+[ycomb,black,thick]
			table[x index=0,y index=1] 
			{g2_j2_r1_beta_12.txt};	\label{symbol:r1beta12_g2}
			\addplot+[ycomb,black,thick]
			table[x index=0,y index=1] 
			{g2_j2_r1_beta_13.txt};	\label{symbol:r1beta13_g2}
			\addplot+[ycomb,black,thick, mark = triangle]
			table[x index=0,y index=1] 
			{g2_j2_r1_beta_23.txt};	\label{symbol:r1beta23_g2}
		\end{axis}
	\end{tikzpicture}
	\caption{Impulse responses $\tilde{g}_{2,1,\beta}(t)$ of $(\tilde{A}_1,b_{2,1,\beta},\tilde{c}_1)$ in Theorem \ref{thm:main2}:  \ref{symbol:r1beta12_g2} $\tilde{g}_{2,1,\{1,2\}}(t)$, \ref{symbol:r1beta13_g2} $\tilde{g}_{2,1,\{1,3\}}(t)$ and \ref{symbol:r1beta23_g2} $\tilde{g}_{2,1,\{2,3\}}(t)$ are strictly positive.}
	\label{fig:SSC_2_impulse_r1_j2}
\end{figure}
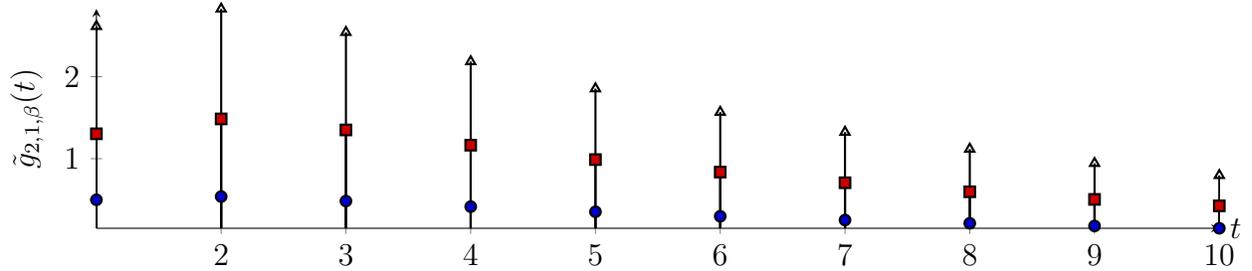

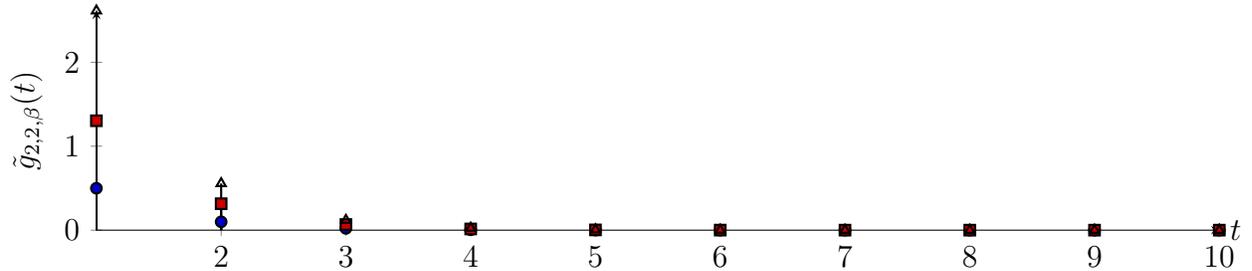
\begin{figure}[t!]
	\centering
	\begin{tikzpicture}
		\begin{axis}[height=4.5cm,
			width=\columnwidth,
			axis y line = left,
			axis x line = center,
			ylabel={$\tilde{g}_{2,2,\beta}(t)$}, 
			xlabel={$t$}, 
			xlabel style={right},
			xmax=10,
			]			
			\addplot+[ycomb,black,thick]
			table[x index=0,y index=1] 
			{g2_j2_r2_beta_12.txt};	\label{symbol:r2beta12_g2}
			\addplot+[ycomb,black,thick]
			table[x index=0,y index=1] 
			{g2_j2_r2_beta_13.txt};	\label{symbol:r2beta13_g2}
			\addplot+[ycomb,black,thick, mark = triangle]
			table[x index=0,y index=1] 
			{g2_j2_r2_beta_23.txt};	\label{symbol:r2beta23_g2}
		\end{axis}
	\end{tikzpicture}
	\caption{Impulse responses $\tilde{g}_{2,2,\beta}(t)$ of $(\tilde{A}_2,b_{2,2,\beta},\tilde{c}_2)$ in Theorem \ref{thm:main2}: \ref{symbol:r2beta12_g2} $\tilde{g}_{2,2,\{1,2\}}(t)$, \ref{symbol:r2beta13_g2} $\tilde{g}_{2,2,\{1,3\}}(t)$ and \ref{symbol:r2beta23_g2} $\tilde{g}_{2,2,\{2,3\}}(t)$ are strictly positive.}
	\label{fig:SSC_2_impulse_r2_j2}
\end{figure}

Next, we consider
\begin{align*}
	A = \begin{pmatrix}
		0.7 &   0.6 &  -2\\
		0.15   & 0.15  & -0.25\\
		0   & 0.03  & 0.1
	\end{pmatrix}, \quad c^\transp = \begin{pmatrix}
		1.1 \\
		0.1\\
		-5.5
	\end{pmatrix}
\end{align*}
with
\begin{equation*}
	\Obs{3} = \begin{pmatrix}
		1.10 &   0.10 &  -5.50\\
		0.79  &  0.51 &  -2.78\\
		0.63   & 0.46 &  -1.98
	\end{pmatrix}
\end{equation*}
Since $\Obs{3}$ contains elements of mixed signs, $\Obs{}$ is not $0$-variation diminishing. Hence, independent of the choice of $b$, no upper bound on the variation of the impulse response can be provided with variation diminishing arguments. Fortunately, as illustrated in Figures~\ref{fig:SSC_2_impulse_r1_j2} \& \ref{fig:SSC_2_impulse_r2_j2}, Theorem \ref{thm:main2} guarantees that $\Obs{}$ is $\svb{1}$. Then, for any $b \in \mathds{R}^3$ with $\vari{b} \leq 1$ we have that $\vari{g} \leq \variz{g} \leq 1$.

\subsection{Example 3: $\svb{1}$ $\Obs{}$ without $k$-positive $A$ representation}\label{subsec:A_not_kpos}
Corollary \ref{cor:k_pos} has the advantage that $A$ does not need to be $k$-positive. This is especially useful in cases where a system does not admit a {realization} with $k$-positive $A$. %
In the following, we derive an example where our results can be used to verify that $\Obs{}$ is $\sr{2}$. However, this system does not permit any realization where $A$ is $\sgc{2}$. Consequently, even a possible extension of Proposition \ref{prop:k_pos_obsv} to the $\sgc{k}$ case would not be applicable to this example. 

We start by defining the system $(\bar{A},\bar{c},\bar{b})$ as
\begin{align*}
	\bar A &:= \begin{pmatrix}
		1  & 0 & 0 & 0 & 0\\
		1  & 1  &  0&0&0\\
		0  & 1 & 1 & 0 & 0\\
            0&0&0& \cos(\textstyle \frac{\pi}{\sqrt{2}}) & -\sin(\textstyle \frac{\pi}{\sqrt{2}}) \\
            0&0&0& \sin(\textstyle \frac{\pi}{\sqrt{2}}) & \cos(\textstyle \frac{\pi}{\sqrt{2}})
	\end{pmatrix}, \quad \bar b := \begin{pmatrix}
		1  \\  1\\ 1 \\ 1 \\ 1
	\end{pmatrix},\\
 \bar c &= \begin{pmatrix}
     1 & 1 & 1 & 0.001 & 0.001
 \end{pmatrix}
\end{align*}
and the corresponding truncated Hankel operator 
\begin{equation*}
    \bar{\mathcal{H}}_{\bar{g}} = \Obs{}(\bar{A},\bar{c}) {\Obs{5}(\bar{A}^\transp,\bar{b}^\transp)}^\transp. %
\end{equation*}
Using the states-space transformation $T = {\Obs{5}(\bar{A}^\transp,\bar{b}^\transp)}^\transp$, we define our pair $(A,c)$ by
\begin{align*}
	A &:= T \bar A T^{-1},\quad c:= \bar c T^{-1}.
\end{align*}
such that
\begin{equation*}
    \Obs{}(A,c) = \bar{\mathcal{H}}_{\bar{g}}.
\end{equation*}
In order to show that $\Obs{}(A,c)$ is $\vd{1}$, it suffices to verify that $\bar{\mathcal{H}}_g$ is $\vd{1}$. To this end, note that 
the impulse response of $(\bar{A},\bar{c}^\transp,\bar{c})$ is given by
\begin{equation*}
  \bar{g}(t) = \frac{t}{2} + 0.002\cos(\frac{\pi}{\sqrt{2}}(t - 1)) + \frac{t^2}{2} + 2,\;  t>0
\end{equation*}
which can be easily verified to be positive. Similarly, the consecutive $2$ minors
\begin{equation*}
     \bar{g}_2(t) = \bar g(t-1) \bar g(t+1) - \bar g(t)^2, \; t > 0
\end{equation*}
can be shown to be negative. The operator resulting from reversing the column order in $\bar{\mathcal{H}}_g$ then has positive $2$-minors and as such fulfills the requirements of Proposition~\ref{prop:consecutive_old}. As $\sgc{2}$ is invariant to column inversion, this shows that $\bar{\mathcal{H}}_g$ is $\sr{2}$ and by Corollary~\ref{cor:sr_k_mat_vd_k_rank_gt_k}, we conclude that $\Obs{}$ is $\vd{1}$. Finally, as the dominant eigenvalues in $\bar{A}$ include complex ones with arguments that are non-rational multiples of $\pi$, it can can be shown by Lemma~\ref{lem:compound_mat} and \cite[Theorem~5]{benvenuti2004tutorial} that there does not exist any realization where $A$ is $\sgc{2}$. 

Note that while for {ease} of exposition, we used Proposition~\ref{prop:consecutive_old} to confirm our claim, it is not difficult to modify this example to cases where $\Obs{}$ is only $\vb{2}$, e.g., by choosing $\bar{c}^\transp = \begin{pmatrix}
    -3 & 1 & 1 & 0.001 & 0.001
\end{pmatrix}$. In such instances, only our new results are applicable, which, however, also come with more involved analytic computations.

\section{Conclusion}
\label{sec:conclusion}
We have derived a tractable approach to certify that the observability operator of a discrete-time LTI system maps a vector with at most $k$ sign changes to a sequence with at most $k$ sign changes. Our derivations are based on new complements to the existing matrix literature of $k$-variation bounding matrices and their application to the observability operator. These complements include the relaxation of rank and dimension assumptions for the characterization of $k$-variation bounding matrices as well as their simplified certification. Interestingly, application of those results to the observability operator transform our problem into a problem of externally positivity verification. While this generalizes recent certificates for the more restrictive notion of order-preserving $k$-variation diminishment \cite{grussler2021internally}, it became also evident that variation bounding and variation diminishing properties in system operators share the same eigenvalue requirements, i.e., the $k+1$ largest eigenvalues of $A$ have to be real and positive.

An advantage of our certificate over \cite{grussler2021internally} is the removal of the $k$-positivity requirement on $A$. In particular, we have constructed an example of a variation-diminishing observability operator, where independent of the system realization such an $A$ does not exist. Our results are applied to the problem of bounding the number of sign changes in the impulse response. 

In future work, we would like to extend our findings to the Hankel and Toeplitz operator without operator splitting and convert our findings into closed-loop design conditions as well as new findings for model order reduction.

\section*{Acknowledgment}
The project was supported by the Israel Science Foundation (grant no. 2406/22), while the first author was also a Jane and Larry Sherman Fellow. 
\appendix
\section*{Appendix}
\section{Proof to Proposition \ref{prop:T_sigma}}
We only prove Items~\ref{item:T_sigma_X_SSC_k} \& \ref{item:T_sigma_X_SVB_k}, as the others follow from \cite[Propositions~5.1.1 \& 5.1.2]{karlin1968total} and \cite[Propositions~2.7]{grussler2021internally}.

\begin{itemize}
    \item Item~\ref{item:T_sigma_X_SSC_k}: from Lemma \ref{lem:compound_mat} Item~\ref{item:Cauchy_Binet}, $T(\sigma)X$ is the product of a positive matrix with a nonnegative/nonpositive matrix. Since all $k$ columns of $X$ are linearly independent, there are no zero columns in $\compound{X}{k}$ and, as such, the product of the matrices results in no zero entries, where all the entries retain the same sign as $\compound{X}{k}$.
    \item Item \ref{item:T_sigma_X_SVB_k}: {From $X$ being $\vb{k-1}$, it holds for any $u$ with $\vari{u}\leq k-1$ that $\vari{Xu} \leq k-1$. Additionally, if $u \neq 0$ then $Xu \neq 0$, due to $\rk(X) = m$. Thus, since $T(\sigma)$ is $\svb{k-1}$ we have $$\variz{T(\sigma) X u} \leq \vari{X u } \leq k-1$$ for any $u\neq 0$ with $\vari{u}\leq k-1$.}
\end{itemize}

\section{Proof to Theorem \ref{thm:main1}}
By Proposition~\ref{prop:n_ssc_mat_svb}, we need to check that $\Obs{s}$ is $\ssgc{n}$, $s \geq n$. Applying Proposition \ref{prop:k_ssc__mat_pena_i} yields that $\Obs{s}$ is $\ssgc{n}$ if and only if the determinants of the following matrices have the same strict sign:
	\begin{equation*}
		M_r[t] = \begin{pmatrix}
			\Obs{n-r} \\
			\Obs{r}A^{n-r + t-1}
		\end{pmatrix}
	\end{equation*}
	for $r\in (1:n) ,\; t \in \mathds{N}_{\geq 1}$.
    For each $r$ the sequence $\{\det(M_r[t])\}_{t \geq 1})$ is the impulse response to $(\tilde{A}_r,\tilde{b}_r,\tilde{c}_r)$ up to the factor of $\det(\left(\Obs{n}\right)^{-1})$. To show this we multiply $M_r[t]$ on the right by $\left(\Obs{n}\right)^{-1}$
    \begin{equation*}
		M_r[t] \left(\Obs{n}\right)^{-1} =\begin{pmatrix}
			\begin{matrix}
				I_{n-r} & 0
			\end{matrix} \\
			\Obs{r} A^{n-r + t -1}\left(\Obs{n}\right)^{-1}
		\end{pmatrix} 
	\end{equation*}
	which is a block-triangular matrix whose determinant can be computed as
	\begin{align*}
		&\det \left(M_r[t]\left(\Obs{n}\right)^{-1}\right) \\
		&=   \det \left(
		\begin{matrix}
			\begin{matrix}
				I_{n-r} & 0
			\end{matrix} \\
			\Obs{r} A^{n -r + t -1}\left(\Obs{n}\right)^{-1}\\
		\end{matrix} 
		\right) \\
		&= \det \left(
		\Obs{r}
		A^{t-1} A^{n-r}
		\left(\Obs{n}\right)^{-1}
		\begin{pmatrix}
			0 \\
			I_r
		\end{pmatrix}
		\right)\\
		&= \pcompound{\Obs{r}}{r}
		\compound{A}{r}^{t-1}
		\pcompound{A^{n-r}\left(\Obs{n}\right)^{-1} \begin{pmatrix} 0 \\ I_r  \end{pmatrix}}{r} 
	\end{align*}
	where the second equality uses the block-diagonal structure and the last equality exploits Cauchy-Binet (see~Lemma \ref{lem:compound_mat}). Then, since 
	\[
	\det \left(M_r[t]\left(\Obs{n}\right)^{-1}\right)
	=
	\det \left(M_r[t]\right)
	\det \left(\left(\Obs{n}\right)^{-1}\right)
	\]
    our claim follows.
	Finally, $\det \left(\Obs{n}\right)$ is a particular minor of $\Obs{}$, which has the same sign as  $\det \left(\left(\Obs{n}\right)^{-1}\right)$, which is why $\det \left(M_r[t]\left(\Obs{n}\right)^{-1}\right)$ needs to be strictly positive.  
\section{Proof to Proposition \ref{prop:ssc_k_mat_svb_k}}
{$\Longleftarrow$:} We begin by noticing that if $X$ is $\svb{k-1}$, then so must be $X_{(1:n,v)}$ for all $v \in \mathcal{I}_{m,j}$ and $j \in (k:m)$. In particular, if $\alpha \in \mathcal{I}_{m,k}$, then by assumption $$\variz{Xu(\alpha)} \leq k-1$$ for all $u(\alpha)$ with $u(\alpha)_i = 0$, $i \not \in \alpha$, {because $$\vari{u(\alpha)} \leq k-1.$$} Hence, by Proposition~\ref{prop:n_ssc_mat_svb}, 
{\begin{equation}
\forall	\alpha \in \mathcal{I}_{m,k}: X_{(1:n),\alpha} \text{ is } \ssgc{k}. \label{eq:ssc_k_all_alpha}
\end{equation}}
In order to conclude that $X$ is $\ssgc{k}$, we need to show strict sign  agreement between all the $X_{(1:n),\alpha}$. {By \eqref{eq:ssc_k_all_alpha}, it suffices to show this only for the first $k$ columns, i.e., that $X_{(1:k),\beta}$ is $\ssgc{k}$ for all $\beta \in \mathcal{I}_{m,k+1}$.} To this end, let $\bar{u} \in \mathds{R}^{k+1}$ be such that $\submatrix{X}{(1:k)}{\beta}\bar{u} = 0$ and $X_{(1:k),\beta_{(1:k)}}\bar{u}_{(1:k)} = X_{(1:k),\beta_{\{k+1\}}}$. Then, by Cramer's rule \cite[Sec. 0.8.3]{horn2012matrix}, 
$$\det(X_{(1:k),\beta_{(1:k)}})\bar{u}_j = (-1)^{k-j}\det(X_{(1:k),\beta^j}), \; $$
where $\beta^j = \beta_{(1:j-1)}\cup\beta_{(j+1:k+1)}$. Our claim is proven if
\begin{equation*}
	\forall j \in (1:k): \det(X_{(1:k),\beta^j})\det(X_{(1:k),\beta^{j+1}}) > 0.
\end{equation*}
Assuming the contrary, there exists at least one $j^\ast$ such that $\bar{u}_{\beta_{j^\ast}}\bar{u}_{\beta{j^\ast + 1}} > 0$. Then, $\vari{\bar{u}} \leq k-1$ and $\variz{X_{(1:n),\beta}\bar{u}} \geq k$, as the first $k$ elements in $X_{(1:n),\beta}\bar{u}$ are zero. This is a contradiction, since $X_{(1:n),\beta}$ is $\svb{k-1}$ by assumption.

{$\Longrightarrow$:} Conversely, let $X$ be $\ssgc{k}$. By Proposition \ref{prop:n_ssc_mat_svb}, $X_{(1:n),\alpha}$ is $\svb{k-1}$ for all $\alpha \in \mathcal{I}_{m,k}$. Thus, for all $u(\alpha)$ as above, $\variz{Xu} \leq k-1$. In order to finish the proof, it remains to consider $u$ with at least $k$ non-zero entries. In this case, $u$ can be partitioned into $k$ parts
	\begin{equation*}
		(u_1,\cdots , u_{s_1}),(u_{s_1 +1},\cdots, u_{s_2}),\cdots ,(u_{s_{k-1}+1},\cdots, u_{m}),
	\end{equation*}
where each part contains no sign changes and at least one non-zero element. Letting $s_0 := 0,\; s_k {:}= m$, we can then define $\tilde{X} \in \mathds{R}^{n \times k}$ and $\tilde{u} \in \mathds{R}^k$ by 
\begin{align*}
\tilde{X}_{(1:n),\{j\}} &:= \sum_{i= s_{j-1} + 1}^{s_j} |u_i| X_{(1:n),\{i\}}, \\
\tilde{u}_j &:= \sign \left(\sum_{i= s_{j-1}+1}^{s_j} u_i \right)
\end{align*}
for $j \in (1:k)$, 
 such that $Xu = \tilde{X} \tilde{u}$. Since the multi-linearity of the determinant implies that also $\tilde{X} \in \mathds{R}^{n \times k}$ is $\ssgc{k}$, it must hold by Proposition \ref{prop:n_ssc_mat_svb} that $$\variz{Xu} = \variz{\tilde{X}\tilde{u}} \leq k-1.$$ 

\section{Proof to Theorem \ref{thm:main2}}
By Proposition~\ref{prop:ssc_k_mat_svb_k}, we need to check that $\Obs{s}$ is $\ssgc{k}$, $s \geq k$. Using Corollary~\ref{cor:k_ssc__mat_pena_ii}, this is equivalent to $\det(O^{N}_{\alpha,\beta})$ sharing the same sign across all defined choices of $\alpha$, $\beta$ and $N$. Since
	\begin{equation*}
	\mathcal{O}^{N}_{\alpha,\beta} = \underbrace{
		\begin{pmatrix}
			\Obs{k-r}  \\
			\Obs{r}A^{k - r + t- 1}
		\end{pmatrix}}_{=: M_{k,r}[t]}
		P_\beta,   
	\end{equation*}
	we need to show that the determinants of $M_{k,r}[t] P_\beta$ coincide with the impulse response of $(\tilde{A}_{r},\tilde{b}_{k,r,\beta},\tilde{c}_r)$. To this end, note that
	\begin{align*}
		&\det \left(M_{k,r}[t] P_\beta \right) = \det \left( M_{k,r}[t]
		\left(\Obs{n}\right)^{-1}\Obs{n}  P_\beta \right)\\
		&=  \det \left(\begin{pmatrix}
			\begin{matrix}
				I_{k-r} & 0
			\end{matrix} \\
			\Obs{r}A^{k - r + t- 1}  \left(\Obs{n}\right)^{-1}
		\end{pmatrix}
		\Obs{n}  P_\beta \right)\\
		&=
		\underbrace{
			\compound{\begin{pmatrix}
					\begin{matrix}
						I_{k-r} & 0
					\end{matrix} \\
					\Obs{r}A^{k - r + t- 1}  \left(\Obs{n}\right)^{-1}
			\end{pmatrix}}{k}}_Q
		\pcompound{ \Obs{n}  P_\beta}{k},
	\end{align*}
	where the last equality uses Lemma~\ref{lem:compound_mat}.
	Further, since $Q$ is block-diagonal, we can express it as
	\begin{align*}
		Q
		&=
		\left(\begin{array}{cc}
			\pcompound{
				\Obs{r}A^{k - r + t- 1}  \left(\Obs{n}\right)^{-1}
				\begin{pmatrix} 0 \\ I_{n-k+r} \end{pmatrix}}{r}
			& 0\end{array}\right) \\
		& = 
		\left(\begin{array}{cc}
			\compound{\Obs{r}}{r}
			\compound{
				A^{t-1}
			}{r}
			\pcompound{
				A^{k-r}\left(\Obs{n}\right)^{-1}
				\begin{pmatrix} 
					0 \\ I_{n-k+r} 
			\end{pmatrix}}
			{r}
			& 0
		\end{array}\right),
	\end{align*}
	where the first equality utilizes the fact that a $k$-minor is zero when the identity matrix $I_{k-r}$ is not part of the corresponding submatrix. 

 Since the considered $\alpha$'s are the same as in Corollary~\ref{cor:k_ssc__mat_pena_ii}, it follows that
	$\Obs{s}$, $s\geq k$, is $\ssgc{k}$ if and only if $(\tilde{A}_{r},\tilde{b}_{k,r,\beta},\tilde{c}_r)$ is strictly externally positive/negative for all $r \in (1:k)$ and all $\beta$ as in Corollary \ref{cor:k_ssc__mat_pena_ii}.

\section{Proof to Corollary \ref{cor:eigen}}
We begin by noting that by Item \ref{item:compound_eigen} in Lemma \ref{lem:compound_mat}, the dominant eigenvalue of $A_{[r]}$ is given by $\prod_{i=1}^r \lambda_i(A)$. Further, since by assumption $\Obs{}$ is $\ssgc{k}$, it follows from Theorem \ref{thm:main2} that for each $\beta \in \mathcal{I}_{n,k}$, the impulse response of the systems defined in Theorem \ref{thm:main2} have to be strictly externally positive/negative. 
Since this requires that each system has a dominant positive pole \cite{ohta1984reachability,farina2011positive,benvenuti2004tutorial}, it follows that $\prod_{i=1}^r \lambda_i(A) \neq 0$ and our claim will be proven by Theorem \ref{thm:main2} if we can show that $\prod_{i=1}^r \lambda_i(A)$ is a controllable and observable mode of $(\tilde{A}_r,\tilde{b}_{k,r,\beta}, \tilde{c}_r)$ for each $r \in (1:k)$ and at least one corresponding choice of $\beta \in \mathcal{I}_{n,k}$.

{Let us first consider the case, when $A$ is diagonalizable. Then, since $(A,c)$ is observable,} there exists an invertible $T \in \mathds{C}^{n \times n}$ with
\begin{align*}
    cT = \begin{pmatrix}
        1 & \cdots & 1
    \end{pmatrix}, \; T^{-1}AT = \begin{pmatrix}
        \lambda_1(A) & & \\
        & \ddots & \\
        & & \lambda_n(A)
    \end{pmatrix}
\end{align*}
and $\lambda_1(A) \neq \lambda_2(A) \neq \dots \neq \lambda_n(A)$. Thus, our claim is proven if $\prod_{i=1}^r \lambda_i(A)$ is an observable and controllable mode of $$(\bar{A}_r,\bar{b}_{k,r,\beta},\bar{c}_r) := (\compound{(T^{-1}AT)}{r},T^{-1}_{[r]} \tilde{b}_{k,r,\beta}, \tilde{c}_r T_{[r]}).$$
Since $\bar{A}_r$ is diagonal with its first element being $\prod_{i=1}^r \lambda_i(A)$, this can only be the case if the first element in $\bar{b}_{k,r,\beta}$ and in $\bar{c}_r$ is non-zero. To see this, observe that since  $\bar{c}_r = (\bar{\mathcal{O}}^r)_{[r]}$ with $\bar{\mathcal{O}}^r = \mathcal{O}^r(\bar{A}_1,\bar{c}_1)$ by Lemma \ref{lem:compound_mat}, it follows from \begin{align*}
   \bar{\mathcal{O}}^r= \begin{pmatrix}
        1 & 1 & \cdots & 1 \\
        \lambda_1(A) & \lambda_2(A) & \cdots & \lambda_n(A)\\
        \vdots & \vdots & & \vdots\\
        \lambda_1^{r-1}(A) & \lambda_2^{r-1}(A) & \cdots & \lambda_n^{r-1}(A) 
    \end{pmatrix}
\end{align*}
and \cite[Sec.~0.9.11]{horn2012matrix} that all elements in $\bar{c}_r$ are non-zero. Further, the first element in $\bar{b}_{k,r,\beta}$ is zero for all $\beta \in \mathcal{I}_{n,k}$ if and only if
\begin{align*}
\underbrace{\begin{pmatrix} \pcompound{L_r\left(\bar{\mathcal{O}}^n \right)^{-1}R_{n-k+r}}{r} & 0 \end{pmatrix} \pcompound{\bar{\mathcal{O}}^n}{k}}_{=:v} \pcompound{ T^{-1}}{k} = 0,
\end{align*}
where $L_j := \begin{pmatrix}
      I_j & 0
  \end{pmatrix} \in \mathds{R}^{j \times n}$ and $R_j : = \begin{pmatrix}
      0 \\
      I_j
      \end{pmatrix} \in \mathds{R}^{n \times j}$. 
Since $\pcompound{ T^{-1}}{k} = T_{[k]}^{-1}$ by Lemma \ref{lem:compound_mat}, this is the case if and only if $v$ = 0. However, by going backwards in the proof of Theorem \ref{thm:main2}, we can see that 
\begin{align*}
    v = \compound{\begin{pmatrix}
					L_{k-r}\bar{\mathcal{O}}^n\\
     L_r \left(\bar{\mathcal{O}}^n \right)^{-1} \bar{\mathcal{O}}^n
			\end{pmatrix}}{k} = \compound{\begin{pmatrix}
					\bar{\mathcal{O}}^n_{(1:k-r),(1:n)}\\
     L_r 
			\end{pmatrix}}{k} \neq 0,
\end{align*}
because $|v_1| = |\det(\bar{\mathcal{O}}^n_{(1:k-r,r+1:k)})| \neq 0$ by \cite[Sec.~0.9.11]{horn2012matrix}. Hence, there exists at least one $\beta$ such that the first element $\bar{b}_r$ is non-zero. {This concludes the proof for diagonalizable $A$.} 
{
Finally, based on the previous case, we will show, now, that in general the eigenvalues $\lambda_1(A),\dots,\lambda_k(A)$ are simple. To this end, assume the opposite, i.e., at least one of the $k$ dominant eigenvalues has algebraic multiplicity greater than 1. Then, by the observability of $(A,c)$, $A$ must be non-diagonalizable (e.g., from the Popov-Belevitch-Hautus test \cite{glad2000control}) such that there exists a $T \in \mathds{C}^{n \times n}$ with
\begin{align*}
	T^{-1}AT = \begin{pmatrix}
		J_{1} & & \\
		& \ddots & \\
		& & J_{l}
		\end{pmatrix} 
\end{align*}
with Jordan blocks of the form 
\begin{align*}
	J_{i} = \begin{pmatrix}
			\lambda_{s_{i-1} + 1}(A) & 1 & & & \\
			& \lambda_{s_{i-1} + 2}(A) & 1 & & \\
			& & \ddots & \ddots & \\
			& & & \lambda_{s_{i}-1}(A) & 1 \\
			& & & & \lambda_{s_i}(A)
		\end{pmatrix}
\end{align*}
with $s_0 = 0$, $s_l = n$, $\lambda_{s_{i-1} + 1}(A)  = \dots = \lambda_{s_i}(A)$ and all $J_i$ pairwise distinct from each other. Next, let a perturbation of $A_{\varepsilon}$ of $A$ be defined by
\begin{align*}
	A_{\varepsilon} := T \begin{pmatrix}
		J_{1,\varepsilon} & & \\
		& \ddots & \\
		& & J_{l,\varepsilon}
	\end{pmatrix} T^{-1} 
\end{align*} 
where if $s_i - s_{i-1} > 1$, then 
\small{\begin{align*}
	J_{i,\varepsilon} := \begin{pmatrix}
		\lambda_{s_{i-1} + 1}(A) & 1 & & & \\
		& \lambda_{s_{i-1} + 2}(A) & 1 & & \\
		& & \ddots & \ddots & \\
		& & & \lambda_{s_{i}-1}(A) & 1 \\
		(-1)^{s_i - s_{i-1}+1} \varepsilon & & & & \lambda_{s_i}(A)
	\end{pmatrix}
\end{align*}}and otherwise $J_{i,\varepsilon} := \lambda_{s_i}(A) -\varepsilon$ for $\varepsilon > 0$. Note that if $\lambda_i(A) \not \in \mathds{R}$, then the corresponding $J_i$ as well as the generalized eigenvector, i.e., the columns of $T$, come in complex conjugates. Therefore, as our perturbation appears in both Jordan blocks at the same position, one can show that $A_{\varepsilon} \in \mathds{R}^{n \times n}$. By Laplace expansion it holds then that 
\begin{equation}
	\det(tI - 	J_{i,\varepsilon}) = (t-\lambda_{s_i}(A))^{s_i-s_{i-1}} + \varepsilon, \label{eq:ew_jordan}
\end{equation}
i.e., each $J_{i,\varepsilon}$ has distinct eigenvalues and as such is diagonalizable. Thus, also $A_\varepsilon$ is diagonalizable and since any minor of the perturbed observability operator $\mathcal{O}(A_{\varepsilon}, c)$ is a continuous function in $\varepsilon$, there exists a $\delta_A > 0$ such that the observability and $\ssgc{k}$ property is retained for all $\varepsilon \in (0, \delta_A)$. By our previous analysis, it follows then that
\begin{equation}
	\lambda_1(A_{\varepsilon}),\dots, \lambda_k(A_{\varepsilon}) \in \mathds{R}_{>0}. \label{eq:real_ew_eps}
\end{equation}
Moreover, if not all $\lambda_1(A),\dots,\lambda_k(A)$ are simple, then there exists a smallest $i^\ast \in \mathds{N}$ satisfying $s_i^\ast-s_{i^\ast-1} > 1$ and $d^\ast := s_{i^\ast-1}+1 \leq k$. Therefore, choosing $\varepsilon \in (0, \min(\delta_A, \argmin_{i\neq j}{\frac{1}{2}|\lambda_i(A) - \lambda_j(A)|^n}))$ implies that $$\{\lambda_j(A_{\varepsilon})\}_{j \in (1:s_{i^\ast-1})} = \{\lambda_j(A)-\varepsilon\}_{j \in (1:s_{i^\ast -1})}$$ and 
\begin{equation*}
	\{\lambda_{s_{i^\ast-1}+j}(A_\varepsilon)\}_{j \in (1:d^\ast)} = \{\lambda_{s_{i^\ast}}(A) + i \varepsilon^{\frac{1}{d^\ast} e^{i \frac{2 \pi j}{d^\ast} }} \}_{j \in (1:d^\ast)}.
\end{equation*}
In particular, latter set contains a real number if and only if $d^\ast$ is odd, which is given by $\lambda_{s_{i^\ast}}(A) -\varepsilon^{\frac{1}{d^\ast}}$. However, if $\lambda_{s_{i^\ast}}(A) -\varepsilon^{\frac{1}{d^\ast}} > 0$, then
\begin{equation*}
	|\lambda_{s_{i^\ast}}(A) -\varepsilon^{\frac{1}{d^\ast}}| < |\lambda_{s_{i^\ast}}(A) + i \varepsilon^{\frac{1}{d^\ast}} e^{i \frac{2\pi}{d^\ast}}|,
\end{equation*}
which contradicts \eqref{eq:real_ew_eps}. 
This shows that $\lambda_1(A),\dots,\lambda_k(A)$ must be simple and by \eqref{eq:real_ew_eps} strictly positive. }

\section{Proof to Theorem \ref{thm:sc_k_mat_vb_k_rank_k}}

Let $X=FG,\; F \in \mathds{R}^{n \times k},\;G \in \mathds{R}^{k\times n}$ be a {rank-revealing} factorization and  $b =Gu \in \mathds{R}^k$. Since $\vari{b}\leq k-1$, we will show the following:
\begin{itemize}
    \item $ii. \implies i.$: if $\vari{Fb}\leq k-1$ then $F$ is $\vb{k-1}$. From Proposition \ref{prop:sc_k_mat_vb_m}, $F$ is then $\sgc{k}$, {or equivalently,} $\compound{F}{k} \in \mathds{R}^{\binom{n}{k}}_{\geq 0} \cup \mathds{R}^{\binom{n}{k}}_{\leq 0}$. Since $\compound{G}{k}^\transp \in \mathds{R}^{\binom{m}{k}}$, it follows by Lemma \ref{lem:compound_mat} Item \ref{item:Cauchy_Binet} that $\compound{X}{k} = \compound{F}{k} \compound{G}{k}$ fulfills our claimed Item \ref{item:column_sgc}.
    
\item $i. \implies ii.$: Since $\rk \left(\compound{X}{k}\right) = 1$ (Lemma \ref{lem:compound_mat} Item \ref{item:rk_k_compound_k}), it follows from {Lemma}~\ref{lem:compound_mat} Item \ref{item:Cauchy_Binet} that $\compound{X}{k}$ is factored as the product of a column vector $\compound{F}{k}$ and row vector $\compound{G}{k}$. By our sign assumption on $\compound{X}{k}$, this requires that $F$ is $\sgc{k}$. 
By Proposition~\ref{prop:sc_k_mat_vb_m} this implies that $F$ is $\vb{k-1}$ and, thus,
\begin{align*}
    \vari{Fb}\leq k-1 & \iff 
    \vari{FGu}\leq k-1\\
    &\iff \vari{Xu} \leq k-1
\end{align*}
which proves our claim. 
\end{itemize}

\section{Proof to Theorem \ref{thm:sc_k_mat_vb_k_rank_gt_k}}
We begin by noticing that $\compound{X}{k}$ has no zero column by assumption,  because any $k$ columns are linearly independent. Further, $T(\sigma)$ is assumed to be as in Lemma~\ref{lem:compound_mat}. Based on that, the following arguments hold. 
\begin{itemize}
    \item $\impliedby$: since  $\pcompound{T(\sigma)X}{k} > 0$ or $\pcompound{T(\sigma)X}{k} < 0$ it follows from Proposition \ref{prop:ssc_k_mat_svb_k} that $T(\sigma)X \; \text{is} \; \svb{k}$, which by Lemma \ref{lem:var_lim} implies that $X \; \text{is} \; \vb{k}$.
\item $\implies$: assuming that $X$ is $\vb{k-1}$, any submatrix, $\submatrix{X}{\nset{n}}{I}$, $I \in \mathcal{I}_{m,k}$ is also $\vb{k-1}$. It follows from Proposition \ref{prop:T_sigma} Item \ref{item:T_sigma_X_SVB_k} then that $T(\sigma)\submatrix{X}{\nset{n}}{I}$ is $\svb{k}$, which is {why} $T(\sigma)\submatrix{X}{\nset{n}}{I}$ is $\ssgc{k}$ by Proposition \ref{prop:ssc_k_mat_svb_k}. Lastly, by Proposition \ref{prop:T_sigma} Item \ref{item:X_SC_k}, $\submatrix{X}{\nset{n}}{I}$ is $\sgc{k}$.\\
\\
It remains to show that each $\pcompound{\submatrix{X}{\nset{n}}{I}}{k}$ consists of elements that share the same sign across all possible choices of $I$. This will be shown by considering a sequence of submatrices with properties that interlink all $\pcompound{\submatrix{X}{\nset{n}}{I}}{k}$. To this end, let $J^* \in \mathcal{I}_{m,k+1}$ be a fixed index set such that $\rk \left(\submatrix{X}{\nset{n}}{J^*}\right) = k+1$ ($J^\ast$ always exists because $\rk(X) > r$). The sequence $\{\submatrix{X}{\nset{n}}{J^i}\}_{i=1}^{k+1}$
of submatrices in $\mathds{R}^{n \times k+1}$ is then defined by a sequence of index sets $\{J^i\}_{i=1}^{k+1}, \; J^i \in \mathcal{I}_{m,k+1}$ with the following properties:
\begin{enumerate}
    \item  $I \subset J^1$ %
    \label{item:first_element}
    \item $J^{k+1} = J^*$.\label{item:last_element}
    \item $\submatrix{X}{\nset{n}}{J^i}$ is full rank for all $i \in \nset{k+1}$
    \item There are $k$ elements in $J^i \cap J^{i+1}$ for all $i \in \nset{k}$ \label{item:overlapping_col}
\end{enumerate}
It is readily seen that such a sequence exists for any $I$:
\begin{equation*}
    J^1 = I \cup j \;\; \text{with} \;\; j\in J^* \; \text{and} \; \rk\left(\submatrix{X}{\nset{n}}{J^1}\right) = k+1
\end{equation*}
Next $J^2$ has one element from $I$ replaced with an element from $J^*$ such that $\rk \left(\submatrix{X}{\nset{n}}{J^2}\right) = k+1$. This is repeated until the sequence is complete. In the case that elements of $I$ overlap with elements of $J^*$, the sequence will be created using the same steps with fewer elements in total.

Now, by construction (Item (\ref{item:overlapping_col})), $\pcompound{\submatrix{X}{\nset{n}}{J^{i}}}{k}$ and $\pcompound{\submatrix{X}{\nset{n}}{J^{i+1}}}{k}$
share an identical column for all $i \in (1:k)$. Since, as argued above, all $\submatrix{X}{\nset{n}}{J^i}$ are $\sgc{k}$, this implies that $\pcompound{\submatrix{X}{\nset{n}}{J^{i}}}{k}$ and $\pcompound{\submatrix{X}{\nset{n}}{J^{i+1}}}{k}$, and as such all $\pcompound{\submatrix{X}{\nset{n}}{J^{i}}}{k}$, $i \in (1:k+1)$, consist of elements that mutually share the same sign. Since  $\submatrix{X}{\nset{n}}{J^\ast}$ is a common element independent of the choice of $I$, our claim follows.

\end{itemize}

\section{Proof to Proposition \ref{prop:row_col_initial_minors_strict_tot_pos}}
Without loss of generality, we will assume that $n \geq m$. We will show that if all row initial column and rows minors of $X$ are positive, then $X$ is strictly totally positive, or equivalently, a submatrix from the first $i$ columns are strictly totally positive by induction. The converse holds trivially. 

The proof is performed by induction over three nested levels:
\begin{enumerate}
    \item on total positivity of submatrices of the initial columns of $X$
    \item on the consecutive $k$-minors including an added added column
    \item on all the consecutive minors of the same order.
\end{enumerate}

\newcounter{inductionCount}

\refstepcounter{inductionCount}
\textbf{Base case \theinductionCount ~($i=1$):} \label{item:Induction1} All elements in $\submatrix{X}{\nset{n}}{\{1\}}$ a part of the column initial minors, which is why $\submatrix{X}{\nset{n}}{\{1\}} > 0$.
 
\textbf{Induction hypotheses \theinductionCount:} Assume that $\submatrix{X}{\nset{n}}{\nset{i-1}}$ is strictly totally positive, $i \in (2:m)$.

\textbf{Induction step \theinductionCount:} We want to show that $\submatrix{X}{\nset{n}}{\nset{i}}$ is strictly totally positive. By Proposition \ref{prop:consecutive_old} and the induction hypotheses, it is sufficient to show that all consecutive minors that include the $i$-th column are strictly positive. This will be done by another induction from the $i$-th order minors down to the first order minors. 

\refstepcounter{inductionCount}
\textbf{Base case \theinductionCount ~($\pcompound{\submatrix{X}{\nset{n}}{\nset{i}}}{i} > 0$):}\label{item:Induction2} This follows from Proposition \ref{prop:col_intial_m_strict_sign_consist} as $X$ is assumed to have positive column initial minors. 

\textbf{Induction hypotheses \theinductionCount:} Let all $i-k+1$ consecutive minors be positive.  

\textbf{Induction step \theinductionCount:} We want to show that all $i-k$ consecutive minors that include the $i$-th column are positive. To see this, we will use another induction over $\det \submatrix{X}{(j:j+i-k-1)}{(k+1:i)}$, $j \in \nset{n-i+k+1}$.

\refstepcounter{inductionCount}
\textbf{Base case \theinductionCount ~($j=1$):}\label{item:Induction3} $\det \submatrix{X}{(1:i-k)}{(k+1:i)}$
is a row initial minor, which by assumption is positive. 

\textbf{Induction hypotheses \theinductionCount:} Assume that the minor $\det \submatrix{X}{(j-1:j+i-k-2)}{(k+1:i)}$ is positive.

\textbf{Induction step \theinductionCount:} Using (\ref{eq:Desnanot_Jacobi_identity}) we have
\begin{align}\label{eq:col_row_init_proof}
    0 < \det \submatrix{X}{(j:j+i-k-1)}{(k+1:i)} = \frac{\alpha_1 \alpha_2 + \alpha_3 \alpha_4}{\alpha_5}
\end{align}
where
\begin{align*}
    \alpha_1 &:= \det \submatrix{X}{(j-1:j+i-k-1)}{(k:i)}\\
    \alpha_2 &:= \det \submatrix{X}{(j:j+i-k-2)}{(k+1:i-1)}\\
    \alpha_3 &:= \det \submatrix{X}{(j-1:j+i-k-2)}{(k+1:i)}\\
    \alpha_4 &:= \det \submatrix{X}{(j:j+i-k-1)}{(k:i-1)}\\
    \alpha_5 &:= \det \submatrix{X}{(j-1:j+i-k-2)}{(k:i-1)},
\end{align*}
because all the minors on the right are positive by induction hypotheses and the assumptions of the theorem. Concretely, induction assumption \ref{item:Induction2} of positive $i-k+1$ consecutive minors gives $\alpha_1$ being positive. From induction assumption \ref{item:Induction1} of strict total positivity of the first $i-1$ columns $\alpha_2, \alpha_4, \alpha_5$ are positive. From induction assumption \ref{item:Induction3}, $\alpha_3$ is positive. This concludes the proof {of} all three inductions and, thus, the claim follows. 

\section{Proof to Proposition \ref{prop:row_col_initial_minors_tot_pos}}
The proof is similar to the proof of Proposition \ref{prop:row_col_initial_minors_strict_tot_pos} with the following changes:
\begin{itemize}
    \item Induction \ref{item:Induction1} is not needed, because we start with $\submatrix{X}{\nset{n}}{\nset{m-1}}$ being strictly totally positive by application of  Proposition \ref{prop:row_col_initial_minors_strict_tot_pos} to the first $m-1$ columns.
    \item From Proposition \ref{prop:consecutive_old} it is sufficient to show that all consecutive minors that are formed by using the $m$-th column need to be positive for minors of orders less than $m$, and nonnegative for minors of order $m$.
    \item For induction \ref{item:Induction2}, the base case still uses Proposition \ref{prop:col_intial_m_strict_sign_consist}, but its non-strict part. 
    \item For induction \ref{item:Induction3}, when showing the sign of \linebreak $\det \submatrix{X}{(j:j+m-k-1)}{(k+1:m)}$, the term $\alpha_1$ is nonnegative instead of positive for $k=1$. 
\end{itemize}

\section{Proof to Theorem \ref{thm:k_sc__mat_pena_i}}
The proof is similar to the proof in \cite[Theorem~2.2(i)]{pena_matrices_1995}. The idea is to use the bijection in Lemma \ref{lem:pena_bijection} between the original matrix $X$ and another matrix, $C$, which has to be totally positive for $X$ to be $\sgc{m}$. Using our modification for the non-strict case in Proposition~\ref{prop:row_col_initial_minors_tot_pos}, this switches a strict sign requirement to a non-strict sign requirement for some of the minors in $C$ and, hence, $X$. %
In particular, the minors in Lemma \ref{lem:pena_bijection} that are allowed to be non-strict are those with 
\begin{align*}
    &r = m\\
    & \gamma_i = m + \alpha_i.
\end{align*}
Since we are looking only at row and column initial minors, all minors are consecutive, which is why $\alpha$ is consecutive. Hence, these correspond to all the consecutive minors of $X$ starting from the $m+1$ row. This concludes the proof.

\section{Proof to Corollary \ref{cor:k_sc__mat_pena_ii}}
{The matrix is $X$ is sliced to consider submatrices $\submatrix{X}{\nset{n}}{\beta}$. For $\beta$'s that are part of the exception case, one applies Theorem \ref{thm:k_sc__mat_pena_i}, and otherwise Proposition \ref{prop:k_ssc__mat_pena_i} can be used.} This then allows applying Theorem \ref{thm:k_sc__mat_pena_i} to all ${\submatrix{X}{\gamma}{\nset{m}}}^\transp$ for $\gamma \in \mathcal{I}_{n,k}$ which gives $\sgc{k}$.

\section{Proof to Proposition \ref{prop:nonstrict_main1}}
The relationship between the minors and the compound systems was already established in the proof of Theorem \ref{thm:main1}. Using Theorem \ref{thm:k_sc__mat_pena_i} we see, now, that the compound system $(\tilde{A}_n,\tilde{b}_n,\tilde{c}_n)$ only needs to be (non-strictly) externally positive, but still requiring that $\tilde{c}\tilde{A}_n^{l-1}\tilde{b} > 0$ for all $l \in (1:n-1)$.

\section{Proof to Proposition \ref{prop:nonstrict_main2}}
From the proof of Theorem \ref{thm:main2} the equivalence between the impulse responses and the minors is already established.
Corollary~\ref{cor:k_sc__mat_pena_ii} can be used then to establish the above relaxed requirement of (strict) external positivity/negativity of the system corresponding the non-strict sign consistency assumption.

\bibliographystyle{plain}

\bibliography{refkpos,refpos,refopt,science}

\end{document}